\documentclass[reqno]{amsart}[standalone]
\usepackage[notocbasic]{nomencl}
\usepackage[table]{xcolor}
\usepackage{tikz}
\usetikzlibrary{positioning, shapes.geometric}
\usepackage{amssymb}
\usepackage{amsmath}
\usepackage{mathrsfs}
\usepackage{amsfonts}
\usepackage{amsmath}
\usepackage{graphicx}
\usepackage{booktabs}
\usepackage{caption}
\usepackage{subcaption}
\usepackage{mathbbol}
\usepackage{hyperref}
\usepackage{lipsum}
\usepackage{mathtools}
\usepackage{multirow}
\usepackage{comment}
\usepackage[english]{babel}
\usepackage{tikz}
\usepackage[algo2e,ruled]{algorithm2e}
\usepackage{algpseudocode}
\graphicspath{{Images/}}
\usepackage[top=1in, bottom=1.25in, left=1.25in, right=1.25in]{geometry}
\usepackage [autostyle, english = american]{csquotes}
\MakeOuterQuote{"}

\usetikzlibrary{babel}
\usepackage[foot]{amsaddr}
\usepackage{pifont}

\usepackage[table]{xcolor}
\usepackage{tikz}
\usepackage{multirow}
\usetikzlibrary{positioning, shapes.geometric}
\usepackage{amssymb}
\usepackage{amsmath}
\usepackage{mathrsfs}
\usepackage[foot]{amsaddr}
\usepackage{amsfonts}
\usepackage{bm}
\usepackage{amsmath}
\usepackage{graphicx}
\usepackage{booktabs}
\usepackage{caption}
\usepackage{mathbbol}
\usepackage{hyperref}
\usepackage{lipsum}
\usepackage{mathtools}

\usepackage{hhline}
\usepackage{pgf}

\DeclareMathOperator*{\argmin}{arg\,min}

 \usepackage{todonotes}

\def\cl {\nonumber \\}
\def\el {\nonumber }

\usepackage{siunitx}

\def\u{{\bm u}}

\def\0{\boldsymbol{0}}

\def\cl {\nonumber \\}
\def\el {\nonumber }

\newtheorem{rem}{Remark}[section]

\begin{document}
\title[OT-Interpolation ROM]{Optimal Transport-Based Displacement Interpolation with Data Augmentation for Reduced Order Modeling of Nonlinear Dynamical Systems}
\author{Moaad Khamlich$^{1}$}
\author{Federico Pichi$^{1}$}
\author{Michele Girfoglio$^{1}$}
\author{Annalisa Quaini$^{2}$}
\author{Gianluigi Rozza$^{1}$}
\address{$^1$ mathLab, Mathematics Area, SISSA, via Bonomea 265, I-34136 Trieste, Italy}
\address{$^2$ Department of Mathematics, University of Houston, 3551 Cullen Blvd, Houston TX 77204, USA}
\begin{abstract}
We present a novel reduced-order Model (ROM)
that leverages optimal transport (OT) theory
and displacement interpolation to enhance
the representation of nonlinear dynamics in complex systems.
While traditional ROM techniques face challenges in this scenario,
especially when data (i.e., observational snapshots) is limited, our method addresses these issues
by introducing a data augmentation strategy based on OT principles.
The proposed framework generates interpolated solutions tracing geodesic paths
in the space of probability distributions, enriching the training dataset
for the ROM.
A key feature of our approach is its ability to provide a continuous
representation of the solution's dynamics
by exploiting a virtual-to-real time mapping.
This enables the reconstruction of solutions at finer temporal scales than those provided by the original data.
To further improve prediction accuracy, we employ Gaussian Process Regression to learn the residual and correct the representation between the interpolated snapshots and the physical solution.

We demonstrate the effectiveness of our methodology with atmospheric mesoscale benchmarks characterized by highly nonlinear, advection-dominated dynamics. Our results show improved accuracy and efficiency in predicting complex system behaviors, indicating the potential of this approach for a wide range of applications in computational physics and engineering.

\end{abstract}
 \maketitle
\section{Introduction}
\label{sec:intro}
The modeling and simulation of complex dynamical systems, which are often characterized by nonlinear dynamics and high-dimensional state spaces, pose significant challenges. Traditional full-order models (FOMs), such as Finite Element and Finite Volume (FV) methods, offer high accuracy but come with significant computational costs.
While high-performance computing facilities can alleviate some of these computational demands \cite{hpc-book}, computing power alone is often not sufficient for real-time simulation, large-scale parametric studies, and efficient uncertainty quantification. Since complex dynamical systems are ubiquitous in Science and Engineering, coupled with tight deadlines and the need to quantify uncertainties in the results, there is an urgent need for more efficient modeling approaches.

Reduced-order models (ROMs) have emerged as a powerful surrogate strategy to alleviate computational bottlenecks, providing rapid evaluations while maintaining accuracy \cite{Benner2015,benner2017model,Grepl2007,schilders,Rozza2008,Lassila2014,Hesthaven2016,Quarteroni2016}. ROMs have been successfully applied across various domains \cite{vol3}, including fluid dynamics \cite{aromabook}, structural mechanics \cite{Farhat2014}, and electromagnetic systems \cite{Wittig2006}. Despite their widespread success, standard ROM techniques face significant challenges when dealing with highly nonlinear, advection-dominated phenomena \cite{Greif2019}. These challenges are particularly evident in systems characterized by moving discontinuities, steep gradients, and traveling waves. Such features often result in a slow decay of the Kolmogorov n-width, which describes the error arising from a projection onto the best-possible linear subspace of a given dimension \cite{Cohen15}.

This slow decay limits the efficacy of linear ROMs and has spurred several research directions, including shifted Proper Orthogonal Decomposition (POD) \cite{Reiss2018,Papapicco2022}, calibration methodologies \cite{nonino2024calibrationbased}, and registration methods \cite{Taddei2020,Ferrero2022}. Additionally, the development of nonlinear ROMs has marked a significant advancement, particularly through the use of deep learning frameworks that efficiently extract relevant features from high-dimensional data \cite{Goodfellow2016}. Novel architectures, such as recurrent neural networks \cite{Maulik2021}, convolutional autoencoders \cite{LeeModelReductionDynamical2020,romor2023,fresca2021} and graph-based approaches \cite{PichiGraphConvolutionalAutoencoder2023,romor2023friedrichs,MorrisonGFNGraphFeedforward2024}, have further refined the potential of neural networks for regression tasks, enabling more flexible and efficient ROMs.

In this paper, we introduce a novel data-driven ROM  for highly nonlinear, advection-dominated systems by incorporating concepts from Optimal Transport (OT) theory \cite{Villani2008}, specifically leveraging displacement interpolation. The ROM is designed to handle both nonlinear dynamics and limited data scenarios, which are common challenges in many practical applications. The key elements of our methodology include an OT-based data augmentation strategy, a time-interpolation parameter mapping, and a correction step to de-bias the barycentric representation. The resulting framework generates physically consistent synthetic data used for a non-intrusive ROM approach able to capture complex nonlinear dynamics efficiently across continuous time scales.

The first step in our methodology is the selection of checkpoints from the original solution snapshots and the computation of OT plans between consecutive checkpoints. These OT plans are crucial for leveraging displacement interpolation to generate synthetic snapshots, allowing for continuous and geometrically coherent transitions between computed solutions. This process enriches the model space with synthetic solutions that are not restricted to the support of the original data snapshots. Moreover, this strategy significantly increases the amount of data available for training, enriching the POD basis, enhancing its expressivity, and improving accuracy at finer temporal scales than those of the original snapshots. Similarly, our approach could be exploited for deep learning-based approaches in model order reduction, where large datasets are crucial but only a limited amount of full-order solutions can be computed.

To construct a ROM capable of continuous inference, we need to establish mappings between the interpolation parameter and physical time. We propose both a linear mapping and an optimal regression approach, balancing flexibility and accuracy. This allows for efficient online prediction at arbitrary time points. As a correction step, we perform POD with regression (POD-R) of the residuals. A Gaussian Process Regression (GPR) model \cite{Rasmussen2004} is trained to predict the POD coefficients of these residuals, compensating for systematic biases in the OT-based interpolation.

The integration of optimal transport with ROMs builds upon recent works in this direction, including interpolation in the general L2-Wasserstein space \cite{mula2020,mula2022}, solution approximation using Wasserstein Barycenters \cite{mula2023,agueh2011barycenters}, alignment of features using OT mappings \cite{Blickhan23}, and the use of the Wasserstein kernel within kPOD reduction \cite{Khamlich23}. Notably relevant to the design of our methodology is work on convex displacement interpolation \cite{Iollo22}, which has been extended in \cite{cucchiara2023model} where nonlinear interpolation is used for data augmentation, focusing on OT maps between Gaussian models. In our approach, we consider the OT problem in its entropic regularized form, popularized through the use of the Sinkhorn algorithm \cite{cuturi13}. Specifically, we use the multiscale formulation proposed in \cite{Schmitzer2019}, which is particularly useful for problems involving vanishing entropy since it overcomes numerical instabilities, slow convergence, and excessive blurring typical of the standard implementations \cite{Peyre2019}.

While the methodology is general and applicable to a wide range of nonlinear dynamical systems, we choose to demonstrate its effectiveness for the simulation of mesoscale atmospheric flow, modeled by the weakly compressible compressible Euler equations. We focus on two classical benchmarks: the smooth rising thermal bubble and the density current. Both benchmarks are challenging as the dynamics is highly nonlinear and there is strong advection dominance, and hence are ideal to test our methodology.

The rest of the paper is organized as follows: Section \ref{sec:rom-model} introduces the general ROM framework of time-dependent problems. Section \ref{sec:ot} provides background on optimal transport and displacement interpolation. Section \ref{sec:ot_rom_framework} details our OT-based ROM framework, including the data augmentation strategy and the non-intrusive nature of the approach. Section \ref{sec:FOM} describes the FOM for the above-mentioned atmospheric benchmarks. Section \ref{sec:numerical_results} presents the numerical results obtained with our ROM and discusses the effectiveness of the methodology. Finally, Section \ref{sec:conclusions} offers conclusions and perspectives for future work.

\section{A Reduced-Order Modeling Framework for Time-Dependent Problems}
\label{sec:rom-model}

Let us start by considering an abstract PDE problem defined over the domain $\Omega \subset \mathbb{R}^n$ and over time interval of interest $\mathcal{I} = [0, t_f] \subset \mathbb{R}$:
\begin{equation}
\label{eq:abs}
\mathcal{L}(u(t)) = 0, \quad t \in \mathcal{I},
\end{equation}
where $u(t) \in \mathcal{V}$ represents the unknown function within an appropriate function space. The operator $\mathcal{L}$ encompasses the differential operators, boundary conditions, and forcing terms that govern the physical behavior of $u(t)$.
For the numerical solution of \eqref{eq:abs}, classical methods are finite difference, finite element, or finite volume methods.
Once discretized in space, the solution at a specific time instant is denoted by the vector $\boldsymbol{u}_h(t) \in \mathbb{R}^{N_h}$, where $N_h$ represents the number of degrees of freedom.

Let $\mathcal{S} = \left\{\boldsymbol{u}_h(t) \mid t \in \mathcal{I} \right\}$ represent the continuous evolution of $\boldsymbol{u}_h$ as time varies. Let us consider a uniform temporal discretization with $N_T$ time steps. We define the time step size as $\Delta t = \frac{t_f}{N_T}$ and the discrete time points as $t_k = k\Delta t$ for $0 \leq k \leq N_T - 1$. With the introduction of time discretization, we obtain a discretized trajectory set as $
    \mathcal{S}_{h} = \{\boldsymbol{u}_{h}^{k}=\boldsymbol{u}_{h}(t_{k})\}_{k=0}^{N_{T}-1} \subset \mathbb{R}^{N_{h}}$.
Our aim is to study $\mathcal{S}_{h}$ to obtain the best possible approximation of $\mathcal{S}$. Note that we can represent $\mathcal{S}_{h}$ in the form of a \textit{snapshot matrix}:
\begin{equation}\label{eq:snapshots_matrix}
\mathbf{S}=\left[\boldsymbol{u}_{h}^{0}\mid\boldsymbol{u}_{h}^{1}\mid\ldots\mid\boldsymbol{u}_{h}^{N_{T}-1}\right] \in \mathbb{R}^{N_{h} \times N_{T}}.
\end{equation}

To reduce the computational cost associated with fine time discretizations, particularly in applications requiring high-frequency monitoring, we are interested in the development of ROMs that leverage regression and interpolation approaches to approximate the solution \textit{continuously} even when using a large time step $\Delta t$.
We begin by recalling the Proper Orthogonal Decomposition, a linear ROM strategy widely used in ROMs. This approach will be exploited both as a metric to compute performance and, together with a regression strategy, to learn residuals.

\subsection{Proper Orthogonal Decomposition (POD)}
\label{sub:pod}

The POD algorithm is used to build a reduced space $U_r$, where $\operatorname{dim}(U_r) = N_r \ll N_h$, spanned by the POD modes, i.e., the columns of the matrix $\mathbf{U_r} \in \mathbb{R}^{N_h \times N_r}$ such that:
\begin{equation}
\mathbf{U_r}=\underset{\hat{\mathbf{U}}\in \mathbb{R}^{N_h \times N_r}}{\operatorname{argmin}} \frac{1}{\sqrt{N_{T}}}\left\|\mathbf{S}-\hat{\mathbf{U}} \hat{\mathbf{U}}^{\top} \mathbf{S}\right\|_{F},
\end{equation}
where $ \|\cdot\|_{F}$ denotes the Frobenius matrix norm and $\mathbf{S}$ is the snapshot matrix in \eqref{eq:snapshots_matrix}.

The POD modes can be computed efficiently through the singular value decomposition (SVD) of the snapshot matrix:
\begin{equation}
\mathbf{S} = \mathbf{U}\mathbf{\Sigma}\mathbf{V}^T,
\end{equation}
where $\mathbf{U} \in \mathbb{R}^{N_h \times N_h}$ and $\mathbf{V} \in \mathbb{R}^{N_T \times N_T}$ are orthogonal matrices, and $\mathbf{\Sigma} \in \mathbb{R}^{N_h \times N_T}$ is a diagonal matrix containing the singular values $\sigma_1 \geq \sigma_2 \geq \cdots \geq \sigma_r > 0$ in descending order. Here, $r$ is the rank of $\mathbf{S}$.

The POD modes are the left singular vectors (i.e., the columns of $\mathbf{U}$) corresponding to the largest singular values. The  eigenvalue energy captured by the first $k$ modes is given by:
\begin{equation}
E(k) = \frac{\sum_{i=1}^k \sigma_i^2}{\sum_{i=1}^r \sigma_i^2}.
\end{equation}
Typically, a reduced basis $\mathbf{U_r} \in \mathbb{R}^{N_h \times N_r}$ is constructed by selecting the first $N_r$ columns of $\mathbf{U}$, where $N_r$ is chosen such that $E(N_r)$ exceeds a desired threshold (e.g., 99.9\% of the total eigenvalue energy).

\subsection{POD with Regression (POD-R)}
\label{sub:pod_r}

The Proper Orthogonal Decomposition with Regression (POD-R) is a data-driven reduction strategy that performs a regression task to approximate the trajectory of the full-order model. This modification not only broadens its utility to encompass scenarios where the FOM solver is not open-source, but also extends its applicability to datasets derived from experimental techniques, such as those obtained from Particle Image Velocimetry (PIV) images \cite{SemeraroAnalysisTimeresolvedPIV2012}.

Using the POD basis $\mathbf{U_r}$, we can approximate the solution $\boldsymbol{u}_h(t)$ by a reduced expansion:

\begin{equation}
	\label{exp-rom}
	\boldsymbol{u}_h(t)\approx \boldsymbol{u}_r(t) =\mathbf{U_r} \boldsymbol{\alpha}_r(t),
\end{equation}
where $\boldsymbol{\alpha}_r \in \mathbb{R}^{N_r}$ is the time-dependent vector collecting the expansion coefficients.
The POD-R method exploits a regression approach based on this linear expansion. Starting from a training database $\mathcal{D} = \{(t_{j}, \boldsymbol{u}_h(t_{j}))\}_{j=1}^{N_\text{train}}$, which consists of $N_\text{train}$ pairs of time points and their corresponding solution, each trajectory is projected onto the reduced basis to obtain the ROM coefficients:

\begin{equation}
	\boldsymbol{\alpha}_{r}(t_j) = \mathbf{U_r}^{\top} \boldsymbol{u}_{h}(t_j), \qquad \forall\ j=1,\dots,N_\text{train}.
\end{equation}
This gives us a database of reduced coefficients, $\mathcal{D}_{r} = \{(t_{j}, \boldsymbol{\alpha}_{r}(t_{j}))\}_{j=1}^{N_\text{train}}$, which can be used to learn a regression map $\pi: \mathcal{I} \mapsto \mathbb{R}^{N_{r}}$, approximating the mapping $\mathscr{F}: t \in \mathcal{I} \rightarrow \boldsymbol{\alpha}_{r}(t) \in \mathbb{R}^{N_{r}}$.

Finally, the full-order solution is reconstructed via the ``inverse'' coordinate transformation:

\begin{equation}
	\label{reconstruction}
	\boldsymbol{u}_h(t) = \mathbf{U_r} \pi(t).
\end{equation}

This method has been used in combination with different regression approaches, including Radial Basis Function (RBF) regression \cite{xiao2015}, Gaussian Process Regression (GPR) and Artificial Neural Networks (ANNs) \cite{ubbiali,PichiArtificialNeuralNetwork2023,HirschNeuralEmpiricalInterpolation2024}. In this study, we exploit the GPR, detailed further in \cite{GuoReducedOrderModeling2018,OrtaliGaussianProcessRegression2022}.
\section{Optimal Transport}
\label{sec:ot}

In the context of the ROM framework introduced thus far, we now turn our attention to Optimal Transport (OT) theory as a powerful tool for enhancing our approach. OT provides a mathematical foundation for interpolating between probability distributions, which we will leverage to generate synthetic snapshots and improve the representation of nonlinear dynamics in our ROM. By incorporating OT-based techniques, particularly displacement interpolation, we aim to address the challenges of capturing complex, time-dependent phenomena and enriching our dataset in scenarios where observational data may be limited. The following section introduces the key concepts of OT that form the basis of our enhanced ROM methodology.

\subsection{Foundations of Optimal Transport Theory}
OT theory \cite{Villani2008,Santambrogio2015} provides a mathematical framework to determine the most cost-effective way to transport one probability distribution into another. Below, we briefly explain that OT, in its discrete form, can be approached as a linear programming problem.

Let $\mu$ and $\nu$ be discrete measures supported on the sets $\mathcal{X} = \{x_{i}\}_{i=1}^{n}$ and $\mathcal{Y} = \{y_{j}\}_{j=1}^{m}$ and defined as follows:
\begin{equation*}
    \mu = \sum_{i=1}^n {a}_i \delta_{x_i}, \quad
    \nu = \sum_{j=1}^m {b}_j \delta_{y_j}, \quad
    \text{with} \quad \boldsymbol{a} \in \mathbb{R}_+^{n}, \quad
    \boldsymbol{b} \in \mathbb{R}_+^{m} \quad \text{s.t.} \quad
    \boldsymbol{a}^{T} \mathbf{1}_n = 1, \quad
    \boldsymbol{b}^{T} \mathbf{1}_m = 1,
\end{equation*}
where $\delta_{x_i}$ and $\delta_{y_j}$ denote Dirac delta functions at points $x_i$ and $y_j$, vectors $\boldsymbol{a}$ and $\boldsymbol{b}$ have non-negative entries representing the weights at $\{x_{i}\}_{i=1}^{n}$ and $\{y_{j}\}_{j=1}^{m}$, and $\boldsymbol{1}_{n}= [1, \ldots, 1] \in \mathbb{R}_{+}^{n}$ denotes the uniform histogram.
We define a cost function $c: \mathcal{X} \times \mathcal{Y} \rightarrow \mathbb{R}^+$, quantifying the unitary transportation cost between elements in sets $\mathcal{X}$ and $\mathcal{Y}$. This function is encapsulated in the cost matrix $\mathbf{C} \in \mathbb{R}_{+}^{n \times m}$, whose entries ${C}_{ij}$  correspond to the cost of transporting mass from $x_i$ to $y_{j}$, raised to the $p$-th power: $C_{ij} = c(x_i, y_j)^p$. The choice of $p$ depends on the specific OT problem formulation, typically $p=1$ or $p=2$.

A \textit{transport plan}\footnote{More generally, a \textit{transport plan} is a probability measure on $\mathcal{X} \times \mathcal{Y}$, whose marginals are $\mu$ and $\nu$. In discrete settings, this is typically represented by a matrix.} or \textit{coupling matrix} $\mathbf{P}$ is a matrix whose entries $P_{ij}$ encode the amount of mass transported from $x_i$ to $y_j$.
A valid transport plan must satisfy the mass conservation constraints, ensuring that the total mass sent from each $x_i$ matches the corresponding weight $a_i$, and the total mass received at each $y_j$ matches $b_j$.
The set of all feasible transport plans is given by:
\[
    \Pi(\mu, \nu) = \left\{ \mathbf{P} \in \mathbb{R}_+^{n \times m} \mid \mathbf{P} \mathbf{1}_m = \boldsymbol{a}, \ \mathbf{P}^T \mathbf{1}_n = \boldsymbol{b} \right\},
\]
and the optimal coupling $\boldsymbol{\pi} \in \Pi(\mu, \nu)$
is defined as the solution that minimizes the total transport cost:
\begin{equation}
\label{eq:unreg}
\boldsymbol{\pi} = \argmin_{\mathbf{P} \in \Pi(\mu, \nu)} \left\langle \mathbf{P}, \mathbf{C} \right\rangle,
\end{equation}
where $\langle \cdot, \cdot \rangle$ denotes the Frobenius inner product between matrices.

Problem \eqref{eq:unreg} can be solved using classical linear programming techniques \cite{transportation-simplex}. However, the computational cost of this approach can quickly become prohibitive as the size of the problems increases \cite{ot-cost}. To solve \eqref{eq:unreg} in high-dimensional settings, one can resort to entropic regularization \cite{cuturi13}, which introduces a smoothing term to the transport problem:
\begin{equation}
\boldsymbol{\pi}_{\epsilon}(\mu, \nu) = \argmin_{\mathbf{P} \in \Pi(\mu, \nu)} [\langle \mathbf{C}, \mathbf{P} \rangle - \epsilon H(\mathbf{P})] \quad \text{with} \quad H(\mathbf{P}) = -\sum_{i=1}^n \sum_{j=1}^m \mathbf{P}_{ij} (\log(\mathbf{P}_{ij})-1),
\end{equation}
where the hyperparameter $\epsilon > 0$ modulates the entropy level, making the plan more diffusive as $\epsilon$ increases. This modification leads to a differentiable and convex optimization problem that can be efficiently solved using the Sinkhorn algorithm \cite{sink}. Since this approach only involves matrix-vector multiplications, it is highly parallelizable and can take advantage of GPUs, resulting in a significant reduction of the computational cost.

\subsection{Displacement Interpolation}
\label{sub:displacement_interpolation}

Let us consider the challenge of defining interpolation between probability measures. An initial approach might be to use a convex combination of two measures $\mu_0$ and $\mu_1$, defined as:
\[
\mu_\alpha = (1-\alpha) \mu_0 + \alpha \mu_1, \quad \alpha \in [0,1] .
\]
However, this linear interpolation is counterintuitive in the context of OT problems because it simply redistributes the mass across the supports of measures $\mu_0$ and $\mu_1$ based on the interpolation parameter $\alpha$, rather than actually transporting mass from one measure to another. Ideally, we would prefer the mass to transition smoothly, effectively moving between the supports.

In the scenario involving Dirac delta measures, a desirable interpolation is given by their McCann interpolation \cite{mccan-interpolation}:
\begin{equation}
    \label{eq:di-deltas}
\mu_\alpha = \delta_{(1-\alpha)x_0 + \alpha x_1} \quad \alpha \in [0,1];
\end{equation}
where $x_0$ and $x_1$ are the locations of the original Dirac deltas.

\begin{figure}[ht]
    \centering
    \includegraphics[width=1\textwidth]{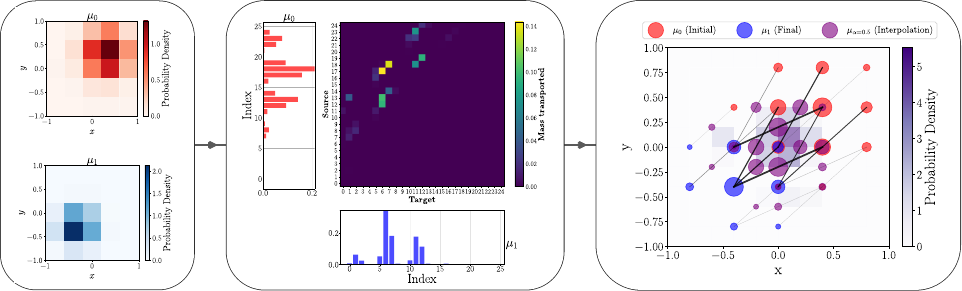}
    \caption{Displacement interpolation between two discrete distributions $\mu_0$ and $\mu_1$ shown on the left panel. In the right panel, the point's masses are proportional to the size of the point and OT plan line thickness are proportional to the mass transported.}
    \label{fig:optimal_transport_di}
\end{figure}

To extend this approach to discrete distributions, consider two discrete probability measures $\mu_0 = \sum_{i=1}^n a_i \delta_{x_i}$ and $\mu_1 = \sum_{j=1}^m b_j \delta_{y_j}$ supported on the same metric space $(\mathcal{X}, d)$. Let $\boldsymbol{\pi} \in \Pi(\mu_0, \mu_1)$ be the OT plan between $\mu_0$ and $\mu_1$, obtained by solving \eqref{eq:unreg}. Then, a generalization of \eqref{eq:di-deltas} is as follows:
\begin{equation}
\label{eq:disp-int}
\mu_\alpha = \sum_{i,j} \pi_{i,j} \delta_{(1-\alpha)x_i + \alpha y_j}  \quad \alpha \in [0,1].
\end{equation}
This interpolation method, known as \textit{displacement interpolation}\footnote{In the most general continuous case, the displacement interpolation is $\mu_\alpha = (e_\alpha)_\# \pi$ where $e_\alpha: \mathcal{X} \times \mathcal{X} \to \mathcal{X}$ is the evaluation map $e_t(x,y) = (1-\alpha)x + ty$, and $(e_\alpha)_\#$ denotes the push-forward operation.}, plays a fundamental role in the field of optimal transport, particularly with regard to the time-dependent version of OT (\cite[Ch.~5]{villani2021topics}). In the discrete-to-discrete setting, this method has an intuitive geometric interpretation: each particle of mass from $\mu_0$ moves towards $\mu_1$ according to the optimal transport plan, taking straight-line paths to their final destinations. See Figure \ref{fig:optimal_transport_di} for an example. At the specified parameter $\alpha$, $\mu_\alpha$ represents the distribution of all particles at their intermediate positions\footnote{In our framework, $\alpha$ is treated as a virtual time parameter, distinct from the physical time $t$. This allows us to interpolate between distributions using $\alpha$ while maintaining a separate notion of actual time evolution in the system.}.

Key properties of displacement interpolation include: (i) geodesicity, as $\mu_\alpha$ traces a geodesic in the Wasserstein space of probability measures; (ii) mass conservation with $\mu_\alpha(X) = 1$ for all $\alpha$; and (iii) finite support at each $\alpha$, accommodating at most $nm$ support points. To the readers interested in more details on displacement interpolation and its properties, we recommend the seminal work by McCann \cite{mccan-interpolation} and the in-depth discussions in Santambrogio's book \cite{Santambrogio2015}, particularly Chapter 5.
\section{OT-Based Interpolation ROM}
\label{sec:ot_rom_framework}

We now present a comprehensive framework for reduced order modeling leveraging displacement interpolation.
This approach combines the theoretical foundations of OT with practical strategies for data extension and efficient online inference, particularly suited for problems with complex, nonlinear dynamics, such as convection-dominated flows described by the model in Sec.~\ref{sec:FOM}.

Let $\mathcal{S}_h = \{\boldsymbol{u}_{h}^{k}=\boldsymbol{u}_{h}(t_{k})\}_{k=0}^{N_{T}-1} \subset \mathbb{R}^{N_{h}}$ be the discretized trajectory set of the FOM solution, as defined in Sec.~\ref{sec:rom-model}.
Our framework aims to extend this discrete set into a continuous representation using OT-based interpolation.

\begin{figure}
    \centering
    \includegraphics[width=\linewidth]{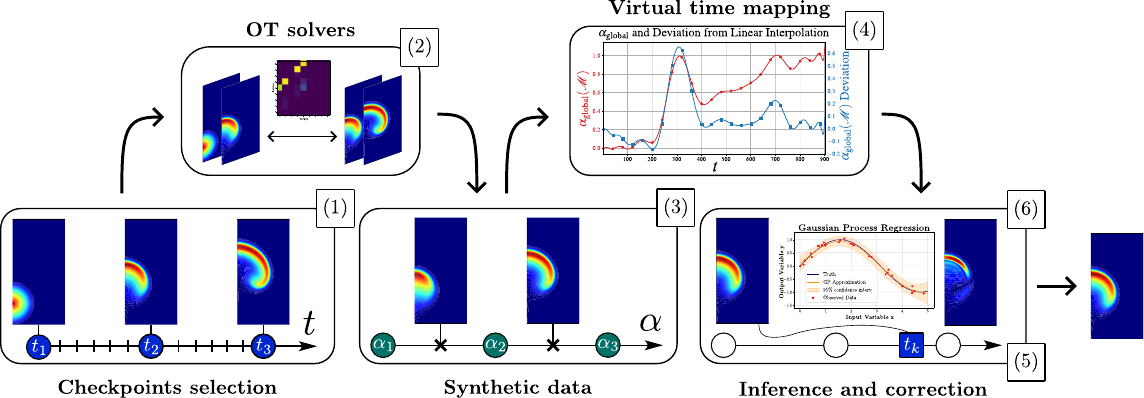}
    \caption{Schematic depiction of our OT-based interpolation ROM approach.}
    \label{fig:scheme}
\end{figure}

Our methodology, depicted in Figure \ref{fig:scheme}, can be subdivided in the following key steps:

\begin{enumerate}
\setlength\itemsep{1em}
\item \textbf{Checkpoint Selection:} we select a subset of $N_c$ checkpoints from $\mathcal{S}_h$, denoted as $\mathcal{C} = \{\boldsymbol{u}_{h}^{k_i}\}_{i=0}^{N_c-1}$, where $0 = k_0 < k_1 < \cdots < k_{N_{c-1}} = N_T -1$. The number of checkpoints $N_c$ is a parameter of the method, allowing for a trade-off between accuracy and computational
efficiency.

\item \textbf{Displacement Interpolation:} for each pair of consecutive checkpoints $(\boldsymbol{u}_{h}^{k_i}, \boldsymbol{u}_{h}^{k_{i+1}})$, we solve the entropic regularized OT problem as described in Sec.~\ref{sub:displacement_interpolation}:
\begin{equation}
    \label{eq:step2}
    \boldsymbol{\pi}_{\epsilon}(\mu_i, \mu_{i+1}) = \argmin_{\mathbf{P} \in \Pi(\mu_i, \mu_{i+1})} [\langle \mathbf{C}, \mathbf{P} \rangle - \epsilon H(\mathbf{P})],
\end{equation}
where $\mu_i$ and $\mu_{i+1}$ are the probability measures corresponding to $\boldsymbol{u}_{h}^{k_i}$ and $\boldsymbol{u}_{h}^{k_{i+1}}$ respectively, and $\epsilon > 0$ is the entropic regularization parameter. For scalar fields, we employ one of two strategies: (i) if the field is non-negative, we normalize it to have unit mass, or (ii) we split the field into positive and negative components, normalizing each separately, and perform the procedure for both.

\item \textbf{Synthetic Data Generation:} Let $\mathcal{X} =\{x_l\}_{l=1}^{N_h}$ be the set of discretization points in our numerical scheme, shared by all fields. These points represent, for example, cell centroids in a FV discretization.
Let $\alpha \in [0,1]$ be the local interpolation parameter, playing the role of virtual time. We define a function $\boldsymbol{u}_{\text{synth}}(i,\alpha)$, for $i \in \{0, 1, \ldots, N_c-2\}$, that generates synthetic snapshots using the OT plan $\boldsymbol{\pi}_{\epsilon}$:
\begin{equation}
    \label{eq:u_synth}
    \boldsymbol{u}_{\text{synth}}(i,\alpha) = \underbrace{\left((1-\alpha)\|\boldsymbol{u}_h^{k_i}\|_1 + \alpha\|\boldsymbol{u}_h^{k_{i+1}}\|_1\right)}_{\text{(A)}} \sum_{l,m} (\boldsymbol{\pi}_{\epsilon})_{lm} \delta_{\text{proj}_{\mathcal{X}}((1-\alpha)x_l + \alpha x_m)},
\end{equation}
where $\|\cdot\|_1$ denotes the L1 norm.
In \eqref{eq:u_synth}, term (A)
interpolates the total mass between consecutive checkpoints and the function $\text{proj}_{\mathcal{X}}(\cdot)$ maps a point to the nearest point in $\mathcal{X}$, ensuring that the interpolated value is associated with the correct discretization point in our numerical scheme\footnote{In our FV implementation, $\text{proj}_{\mathcal{X}}(\cdot)$ specifically maps a point to the centroid of its containing cell. This ensures that the interpolated value is correctly associated with the cell-averaged quantities typical in FV discretizations.}.

We generate a set of $N_{\text{synth}}$ synthetic snapshots within each interval $[t_{k_i}, t_{k_{i+1}}]$:
\begin{equation}
    \boldsymbol{u}_{\text{synth}}^{i,j} = \boldsymbol{u}_{\text{synth}}(i, \alpha_j), \quad \text{where} \quad \alpha_j = j/(N_{\text{synth}}+1) \quad \text{for} \quad j = 1, \ldots, N_{\text{synth}}.
\end{equation}
These synthetic snapshots enrich our dataset, providing a denser sampling of the solution manifold for the time-interpolation parameter mapping and the subsequent POD step, which captures the dominant modes of the system's behavior.

\item \textbf{Time-Interpolation Parameter Mapping:} We seek a mapping between the physical time $t$ and the interpolation parameter $\alpha$, to label the synthetic snapshots with ``real'' instead of ``virtual'' time. To better understand this concept, we define $\alpha_{\text{local}} \in [0,1]$ as the interpolation parameter within a specific interval between two checkpoints, and $\alpha_{\text{global}} \in [0,1]$ as the overall progress through the entire simulation.
For equispaced checkpoints with interval $\Delta t_c = \frac{t_f}{N_c - 1}$, we express $\alpha_{\text{global}}$ as:
\begin{equation}
\alpha_{\text{global}} = \frac{i\Delta t_c}{t_f} + \alpha_{\text{local}}\frac{\Delta t_c}{t_f} = \frac{i + \alpha_{\text{local}}}{N_c - 1}, \quad i \in \{0, 1, ..., N_c - 2\},
\end{equation}
where $i$ is the index of the current checkpoint interval. The first term represents the fraction of the total real time that has elapsed up to the start of the current checkpoint interval, while the second term represents the progress within the current interval.

We aim to find a mapping $\mathcal{F} = (\mathcal{F}_1,\mathcal{F}_2): [0, t_f] \rightarrow \{0,\ldots,N_c-2\} \times [0, 1]$ that returns the interval index $i$ and the local interpolation parameter $\alpha_{\text{local}}$ for any given time $t$. We propose two strategies:
\begin{itemize}
    \item[-] \textbf{Linear Mapping:} We define a piecewise linear mapping $\mathcal{L}: [0, t_f] \rightarrow \{0,\ldots,N_c-2\} \times [0, 1]$ as:
\begin{equation}
    \mathcal{L}(t) = \left(\left\lfloor\frac{t}{\Delta t_c}\right\rfloor,\ \frac{t - t_{k_i}}{t_{k_{i+1}} - t_{k_i}}\right), \quad t \in [t_{k_i}, t_{k_{i+1}}],
\end{equation}
assuming a uniform progression of $\alpha_{\text{local}}$ with respect to time within each checkpoint interval and set $\mathcal{F} = \mathcal{L}$.
\item[-] \textbf{MinL2 Regression:} We find the optimal mapping $\mathcal{M}: [0, t_f] \rightarrow \{0,\ldots,N_c-2\} \times [0, 1]$ by solving:
\begin{equation}
\label{eq:minL2}
    \mathcal{M}(t) = \argmin_{(g_1, g_2)} \sum_{k=0}^{N_T} \|\boldsymbol{u}_h^k - \boldsymbol{u}_{\text{synth}}(g_1(t_k), g_2(t_k))\|_2^2,
\end{equation}
where $g_1: [0, t_f] \rightarrow \{0,\ldots,N_c-2\}$ determines the interval index and $g_2: [0, t_f] \rightarrow [0, 1]$ determines the local interpolation parameter. In practice, we solve problem \eqref{eq:minL2} using a dictionary-based approach.
First, for each discrete time point $t_k$, we compute:
\begin{equation}
\alpha_{\text{global},k} =  \frac{i^* + j^*/(N_{\text{synth}}+1)}{N_c -1}, \quad (i^*, j^*) = \argmin_{(i,j)} ||\boldsymbol{u}_h(t_k) - \boldsymbol{u}_{\text{synth}}^{i,j}||_2^2.
\end{equation}
Then, we perform a regression on the resulting database $\{(t_k, \alpha_{\text{global},k})\}_{k=0}^{N_T-1}$ to obtain a continuous function $\tilde{\alpha}_{\text{global}}: [0, t_f] \rightarrow [0, 1]$.
For any given time $t$, we can compute the interval index $\tilde{i}$ and local interpolation parameter $\tilde{\alpha}_{\text{local}}$ from the predicted $\tilde{\alpha}_{\text{global}}(t)$ as follows:
\begin{equation}
\begin{aligned}
    \tilde{i} &= \lfloor \tilde{\alpha}_{\text{global}}(t) \cdot (N_c - 1) \rfloor, \quad \text{and} \quad   \tilde{\alpha}_{\text{local}} &= \tilde{\alpha}_{\text{global}}(t) \cdot (N_c - 1) - \tilde{i},
\end{aligned}
\end{equation}
expressing the MinL2 regression mapping as $\mathcal{M}(t) = (\tilde{i}, \tilde{\alpha}_{\text{local}})$. Then, we set $\mathcal{F} = \mathcal{M}$.
\end{itemize}

\item \textbf{Online Inference:} For a given time $t^*$, we predict the solution $\boldsymbol{u}_h(t^*)$ using the chosen mapping $\mathcal{F}$:
\begin{equation}
    \label{eq:online_inference}
    \boldsymbol{u}_h(t^*) \approx \boldsymbol{u}_{\text{synth}}(\mathcal{F}_1(t^*), \mathcal{F}_2(t^*)),
\end{equation}
where $\mathcal{F}_1$ and $\mathcal{F}_2$ are the components of $\mathcal{F}$. This step allows for rapid prediction of system states at arbitrary time points, leveraging the pre-computed OT plans and time-$\alpha$ mappings.
\item \textbf{POD-based Correction:} To enhance prediction accuracy, we introduce a correction step based on POD with regression (POD-R) of the residuals. Let $\mathbf{R} = [\boldsymbol{r}_0 | \cdots | \boldsymbol{r}_{N_T-1}] \in \mathbb{R}^{N_h \times N_T}$ be the matrix of residuals, where $\boldsymbol{r}_k = \boldsymbol{u}_h^k - \boldsymbol{u}_{\text{synth}}(\mathcal{F}_1(t_k), \mathcal{F}_2(t_k))$ for $k = 0, \ldots, N_T-1$. We perform POD on $\mathbf{R}$ to obtain a reduced basis $\mathbf{U}_r \in \mathbb{R}^{N_h \times N_r}$. The columns of $\mathbf{U}_r$ form a basis for the dominant modes of the residuals, capturing systematic biases in the OT-based interpolation.

We train a Gaussian Process Regression (GPR) model $\mathcal{G}: [0, t_f] \rightarrow \mathbb{R}^{N_r}$ to predict the POD coefficients of the residuals:
\begin{equation}
    \delta \boldsymbol{\alpha}_r(t) = \mathcal{G}(t),
\end{equation}
where $\delta \boldsymbol{\alpha}_r(t)$ are the POD coefficients of the residual at time $t$. Then, the corrected prediction for $\boldsymbol{u}_h(t^*)$ is given by:
\begin{equation}
    \label{eq:corrected_prediction}
    \boldsymbol{u}_h(t^*) \approx \boldsymbol{u}_{\text{synth}}(\mathcal{F}_1(t^*), \mathcal{F}_2(t^*)) + \mathbf{U}_r \delta \boldsymbol{\alpha}_r(t^*).
\end{equation}
This correction compensates for the bias in the barycentric representation, improving the alignment with the true solution manifold.
\end{enumerate}

\begin{rem}
\label{remark1}
It is worth noting that, for each checkpoint, $\alpha_{\text{local}}$ is either 0 (at the start of an interval) or 1 (at the end of an interval). Therefore, the corresponding $(t,\alpha_{\text{global}})$ pairs at the checkpoints can be expressed as:
\begin{equation}
\left(t_{k_i}, \alpha_{\text{global}}(t_{k_i})\right) = \left(i \Delta t_c, \frac{i}{N_c - 1}\right), \quad 0 \le i \le N_c -1,
\end{equation}
where $i$ represents the checkpoint index. These points form a diagonal line in the $[0,t_f] \times [0,1]$ space.
\end{rem}

\begin{rem}
\label{remark2}
As $N_c$ increases, $\alpha_{\text{global}}(\cdot,\mathcal{M})$ converges towards $\alpha_{\text{global}}(\cdot,\mathcal{L})$. This is due to the increased density of fixed points where both mappings must agree, reducing the scope for $\mathcal{M}$ to deviate from $\mathcal{L}$ between checkpoints. Consequently, local nonlinearities captured by $\mathcal{M}$ become negligible compared to the overall linear trend.
\end{rem}

While the complete framework provides a comprehensive ROM model, the initial steps (1)-(3)
alone provide a \emph{data augmentation} technique
that can be valuable for the training of different
non-intrusive ROMs.
The synthetic data generation process allows us to construct a comprehensive snapshot matrix that captures both the original checkpoints and the interpolated states. We define this synthetic snapshot matrix $\mathbf{S}_{\text{synth}} \in \mathbb{R}^{N_h \times N_\text{tot}}$, where $N_\text{tot} = N_\text{synth}(N_c - 1) + N_c $ is the total number of snapshots including the checkpoints:
\begin{equation}
    \mathbf{S}_{\text{synth}} = [\boldsymbol{u}_{h}^{k_0} | \boldsymbol{u}_{\text{synth}}^{0,1} | \cdots | \boldsymbol{u}_{\text{synth}}^{0,N_\text{synth}} | \boldsymbol{u}_{h}^{k_1} | \cdots | \boldsymbol{u}_{h}^{k_{N_c-1}}].
\end{equation}
This enriched dataset provides a denser sampling of the solution manifold, offering a more comprehensive representation of the system's evolution.

A possible use of $\mathbf{S}_{\text{synth}}$ is as the basis for POD to $\mathbf{S}_{\text{synth}}$ to obtain a reduced basis $\mathbf{U}_r \in \mathbb{R}^{N_h \times N_r}$:
\begin{equation}
    \mathbf{U}_r = \text{POD}(\mathbf{S}_{\text{synth}}, N_r).
\end{equation}
The resulting reduced basis and augmented dataset provide a rich foundation for developing various non-intrusive ROM techniques, potentially improving their accuracy and efficiency by incorporating knowledge from both the original checkpoints and the OT-based interpolated states.
\section{Problem Definition and Its Numerical Discretization}\label{sec:FOM}

If we neglect the effects of moisture, solar radiation, and ground heat flux, the dynamics of the dry atmosphere can be described by the Euler equations:
\begin{align}
&\frac{\partial \rho}{\partial t} + \nabla \cdot (\rho \u) = 0 &&\text{in} \ \Omega \ \times \ (0, t_f], \label{eq:continuity} \\
&\frac{\partial (\rho \u)}{\partial t} + \nabla \cdot (\rho \u \otimes \u) + \nabla p + \rho g \widehat{\mathbf{k}} = \mathbf{0} &&\text{in} \ \Omega \ \times \ (0, t_f], \label{eq:momentu} \\
&\frac{\partial (\rho e)}{\partial t} + \nabla \cdot (\rho \u e) + \nabla \cdot (p \u) = 0 &&\text{in} \ \Omega \ \times \ (0, t_f].\label{eq:energy}
\end{align}
where the unknowns are the air density $\rho$, the wind velocity $\u$, the pressure $p$, and the total energy density $e$, which is defined as:
\begin{equation}
e = c_v T + \frac{1}{2}|\u|^2 + gz, \label{eq:total_energy_density}
\end{equation}
where $T$ is the temperature and $c_v=715.5$ J/(Kg K) is the specific heat capacity at constant volume. Additionally, $g$ is the gravitational constant and $\widehat{\mathbf{k}}$ is the unit vector in the vertical $z$ direction. Eqs.~\eqref{eq:continuity}-\eqref{eq:energy} state conservation of mass, momentum, and total energy, respectively.
We assume that dry air behaves as an ideal gas \cite{Pielke2002}. Then, for closure we add to  system \eqref{eq:continuity}-\eqref{eq:energy} the ideal gas law:
\begin{align}
p = \rho R T \quad \text{in} \ \Omega \ \times \ (0, t_f], \label{eq:EQ_of_State}
\end{align}
where $R= 287$ J/(Kg K) is the specific gas constant for dry air.

For numerical stability, we add an artificial diffusion term to the momentum and energy equations:
\begin{align}
&\frac{\partial (\rho \u)}{\partial t} +  \nabla \cdot (\rho \u \otimes \u) + \nabla p + \rho g \widehat{\mathbf{k}}  - \mu_a \Delta \u = \mathbf{0}  &&\text{in } \Omega \times (0,t_f], \label{eq:momentum_stab}\\
&\frac{\partial (\rho e)}{\partial t} +  \nabla \cdot (\rho \u e) +   \nabla \cdot (p \u) - c_p \dfrac{\mu_a}{Pr} \Delta T = 0 &&\text{in } \Omega \times (0,t_f], \label{eq:energy_stab}
\end{align}
where $\mu_a$ is a constant (artificial) diffusivity coefficient and $Pr$ is the Prandtl number. We note that this is equivalent to using a basic eddy viscosity model in Large Eddy Simulation (LES) \cite{BIL05}. More sophisticated LES models can be found in, e.g., \cite{CGQR, marrasNazarovGiraldo2015}.

It is typical to decompose the pressure into a  baseline state $p_0$, coinciding with the pressure in hydrostatic balance, and a perturbation $p^{\prime}$:
\begin{align}
p = p_0 + p'. \label{eq:EQ_p_splitting}
\end{align}
We can apply a similar decomposition to the density:
\begin{align}
\rho = \rho_0 + \rho'. \label{eq:rhoPrime}
\end{align}
Note that from eq. \eqref{eq:EQ_of_State} we have:
\begin{align}
\rho_0 = \dfrac{p_0}{R T_0}.
\label{eq:rho_hyd}
\end{align}
and the assumption of hydrostatic balance gives:
\begin{align}
\nabla p_0 + \rho_0 g \widehat{\mathbf{k}} = \boldsymbol{0}. %
\label{eq:p0_hydro}
\end{align}
By plugging \eqref{eq:EQ_p_splitting} into \eqref{eq:momentum_stab} and  accounting for \eqref{eq:rho_hyd}, we can recast the momentum conservation equation into:
\begin{align}
&\frac{\partial (\rho \u)}{\partial t} +  \nabla \cdot (\rho \u \otimes \u) + \nabla p' + \rho' g \widehat{\mathbf{k}} - \mu_a \Delta \u = \mathbf{0}  &&\text{in } \Omega \times (0,t_f]. \label{eq:momeCE2}%
\end{align}
Concerning eq.~\eqref{eq:energy_stab}, it is worth noting that a quantity of interest for atmospheric problems is the potential temperature
\begin{equation}
\label{eq:potential_temperature}
\theta = \frac{T}{\pi}, \quad \pi = \left( \frac{p}{p_g} \right)^{\frac{R}{c_p}},
\end{equation}
where $\pi$ is the so-called Exner pressure, $p_g = 10^5$ Pa is standard pressure at the ground, and $c_p = c_v + R$ is the specific heat capacity at constant pressure for dry air.
In order to avoid computing \(\theta\) in a post-processing phase, one can replace eq.~\eqref{eq:energy_stab} with:
\begin{equation}
\label{eq:energy3}
\frac{\partial (\rho \theta)}{\partial t} + \nabla \cdot (\rho \u \theta)  - \dfrac{\mu_a}{Pr} \Delta \theta
= 0 \quad \text{in} \ \Omega \ \times \ (0, t_f].
\end{equation}
Equation of state \eqref{eq:EQ_of_State} in terms of $\theta$ instead of $T$ is:
\begin{equation}
\label{eq:state_potential}
p = p_g \left( \frac{\rho R \theta}{p_g} \right)^{\frac{c_p}{c_v}}.
\end{equation}
The use of potential temperature offers several advantages, especially for adiabatic atmospheric processes. For instance, in adiabatic hydrostatic conditions, the potential temperature remains constant with respect to the height from the ground and equal to the absolute atmospheric temperature at the ground level.

The system to be solved is given by eqs.~\eqref{eq:continuity}, \eqref{eq:EQ_p_splitting}, \eqref{eq:rhoPrime}, \eqref{eq:momeCE2}, \eqref{eq:energy3} and \eqref{eq:state_potential}, which needs to be supplemented with appropriate initial and boundary conditions.
For more details on the formulation of the Euler equations adopted in this work and a comparison of different formulations, we refer the interested reader to, e.g.,  \cite{girfoglio2024comparative}.
In order to efficiently solve the coupled problem, we adopt a splitting approach thoroughly described in \cite{Girfoglio_2023, girfoglio2024comparative, CGQR}, which uses Backward Differencing Formula of first order (BDF1) for the approximation of the time derivative, a semi-implicit scheme for the convective terms,  a fully implicit method for the diffusion terms, and a second-order Finite Volume method for space discretization.

For more details about the numerical method, the reader is referred to \cite{Girfoglio_2023,girfoglio2024comparative,CGQR}.
\section{Numerical Results}
\label{sec:numerical_results}

To validate our framework, we consider two standard benchmarks for atmospheric flows in two dimensions: the rising thermal bubble with the setting from \cite{Ahmad2006} and the density current with the setting from \cite{Straka1993}.
Both test cases model a perturbation of a neutrally stratified atmosphere with uniform background potential temperature. Sec.~\ref{subsec:thermal-bubble} presents the results for the rising thermal bubble benchmark, in which the perturbation consists of warm air that rises and deforms into a mushroom shape.
Sec.~\ref{subsec:density-current} shows the results for the density current benchmark, in which the perturbation is made of cold air that falls to the ground and forms a propagating cold front.

To quantify the accuracy of our approach, we use the relative $\mathcal{L}^2$ norm in spatial domain $\Omega$, defined as:
\begin{equation}
\label{eq:l2_error}
E_{\mathcal{L}^2(\Omega)}(u_{\text{ref}}, \hat{u}) = \frac{\|u_{\text{ref}} - \hat{u}\|_{\mathcal{L}^2(\Omega)}}{\|u_{\text{ref}}\|_{\mathcal{L}^2(\Omega)}},
\end{equation}
where $u_{\text{ref}}$ is the reference solution field and $\hat{u}$ is an approximation. For both tests, the quantity of interest is the perturbation of potential temperature $\theta'$. The reference solution $\theta'_{\text{ref}}(t)$ is obtained with a fine time step $\Delta t_{\text{ref}}$.

To generate the approximation, we introduce a second coarser time step $\Delta t_{\text{train}} = m \Delta t_{\text{ref}}$, for some integer $m > 1$, and denote with $\theta'_{\text{train}}(t)$ the training solution.
We define the set of training time points as $\mathcal{T}_{\text{train}} = \{t_k = k\Delta t_{\text{save}} \mid k = 1, ..., N_T\}$, where $\Delta t_{\text{save}} = n \Delta t_{\text{train}}$ for some integer $n \geq 1$.
The corresponding training dataset is given by $\mathcal{S}_{\text{train}} = \{\theta'_{\text{train}}(t_k) \mid t_k \in \mathcal{T}_{\text{train}}\}$.
To evaluate our method's capability across different time discretizations, we define a set of test points $\mathcal{T}_{\text{test}}$ that has no intersection with $\mathcal{T}_{\text{train}}$ and
denoted with $\mathcal{S}_{\text{test}} = \{\theta'_{\text{ref}}(t_j) \mid t_j \in \mathcal{T}_{\text{test}}\}$
the corresponding test dataset.

The following error metrics will be considered.
The coarser discretization introduces an intrinsic error quantified by the time \emph{discretization error}:
\begin{align}
E_{\text{disc}}(t_k, \Delta t_{\text{train}}) &= E_{\mathcal{L}^2(\Omega)}(\theta'_{\text{train}}(t_k),\theta'_{\text{ref}}(t_k)), \quad t_k \in \mathcal{T}_{\text{train}}, \cl
\overline{E}_{\text{disc}}(\Delta t_{\text{train}}) &= \frac{1}{|\mathcal{T}_{\text{train}}|} \sum_{t_k \in \mathcal{T}_{\text{train}}} E_{\text{disc}}(t_k, \Delta t_{\text{train}}), \el
\end{align}
where $E_{\text{disc}}$ is the pointwise error and $\overline{E}_{\text{disc}}$ is its mean across time.

The \textit{interpolation error} quantifies the accuracy of our synthetic snapshots interpolating the training data:
\begin{equation}
E_{\text{interp}}(\mathcal{F}, t_k) = E_{\mathcal{L}^2(\Omega)}(\boldsymbol{u}_{\text{synth}}(\mathcal{F}_1(t_k), \mathcal{F}_2(t_k)), \theta'_{\text{train}}(t_k)), \quad t_k \in \mathcal{T}_{\text{train}} ,
\end{equation}
where $\mathcal{F} \in \{\mathcal{M}, \mathcal{L}\}$ represents either the MinL2 regression mapping $\mathcal{M}$ or the linear mapping $\mathcal{L}$, as explained in Sec.~\ref{sec:ot_rom_framework}, and
$\boldsymbol{u}_{\text{synth}}$ is the synthetic snapshot generation function defined in eq.~\eqref{eq:online_inference}.
The \textit{generalization error} measures how well our synthetic snapshots approximate the reference solution:
\begin{equation}
E_{\text{gen}}(\mathcal{F}, t_j) = E_{\mathcal{L}^2(\Omega)}(\boldsymbol{u}_{\text{synth}}(\mathcal{F}_1(t_j), \mathcal{F}_2(t_j)), \theta'_{\text{ref}}(t_j)), \quad t_j \in \mathcal{T}_{\text{test}}.
\end{equation}
It is important to note that the generalization error represents a cross-discretization comparison, evaluating our method's ability to approximate the finer reference solution using coarser training data for time instances unseen during training.
For both error types, we report the mean error:
\begin{equation}
\overline{E}_{\text{type}}(\mathcal{F}) = \frac{1}{|\mathcal{T}|} \sum_{t \in \mathcal{T}} E_{\text{type}}(\mathcal{F}, t), \quad \text{type} \in \{\text{interp}, \text{gen}\},
\end{equation}
where $\mathcal{T}$ is the corresponding set of time points for each error type.

\vskip .2cm
\noindent {\bf Computational Setup}.
For the numerical solution of the rising thermal bubble and density current benchmarks, we  used
GEA\footnote{The source code of GEA can be found in \url{https://github.com/GEA-Geophysical-and-Environmental-Apps/GEA}} - Geophysical and Environmental Applications \cite{GEA,girfoglio2024comparative, Girfoglio_2023}, an open-source software package for atmospheric and oceanic modeling based on
the OpenFOAM library \cite{openfoam}.
For the OT problems, we used the MultiScaleOT library \cite{multiscale-ot} which implements stabilized sparse scaling algorithms for entropy regularized transport problems \cite{Schmitzer2019}.
The experiments were performed on a machine equipped with an NVIDIA GeForce RTX 3080 GPU and an 11th Gen Intel(R) Core(TM) i7-11700 @ 2.50GHz processor running Rocky Linux version 8.8.

\subsection{Rising Thermal Bubble}
\label{subsec:thermal-bubble}
\graphicspath{{Images/thermal_bubble/}}

The computational domain for this benchmark is $\Omega = [\SI{0}{\meter}, \SI{5000}{\meter}] \times [\SI{0}{\meter}, \SI{10000}{\meter}]$ in the $xz$-plane. We discretize the domain using a uniform structured mesh with mesh size of $h = \Delta x = \Delta z = \SI{125}{\meter}$, resulting in a $40 \times 80$ grid and $N_h = 3200$.

The initial conditions are defined as follows:
\begin{equation}
\left\{
\begin{aligned}
\theta_s &= \SI{300}{\kelvin}, \\[6pt]
\theta_0 &= \begin{cases}
    \theta_s + 2\left[1 - \frac{r}{r_0}\right] & \text{if } r \leq r_0 = \SI{2000}{\meter} \\
    \theta_s & \text{otherwise}
\end{cases} ,\\[6pt]
\rho_0 &= \frac{p_g}{R\theta_0}\left(\frac{p}{p_g}\right)^{c_v/c_p} \text{ with } p = p_g\left(1 - \frac{gz}{c_p\theta_0}\right)^{c_p/R}, \\[6pt]
\u_0 &= \mathbf{0}, \\[6pt]
\end{aligned}
\right.
\end{equation}
where $r = \sqrt{(x - x_c)^2 + (z - z_c)^2}$ is the distance from the center of the perturbation $(x_c, z_c) = (\SI{0}{\meter}, \SI{2000}{\meter})$. The values of constants $p_g$, $g$, $c_p$, and $R$ are as given in Sec.~\ref{sec:FOM}. This perturbation in the potential temperature field affects the pressure and density fields through the adiabatic balance equation \eqref{eq:potential_temperature} and the equation of state \eqref{eq:state_potential}, respectively.

\begin{figure}[ht]
	\centering
    \includegraphics[width=.8\textwidth]{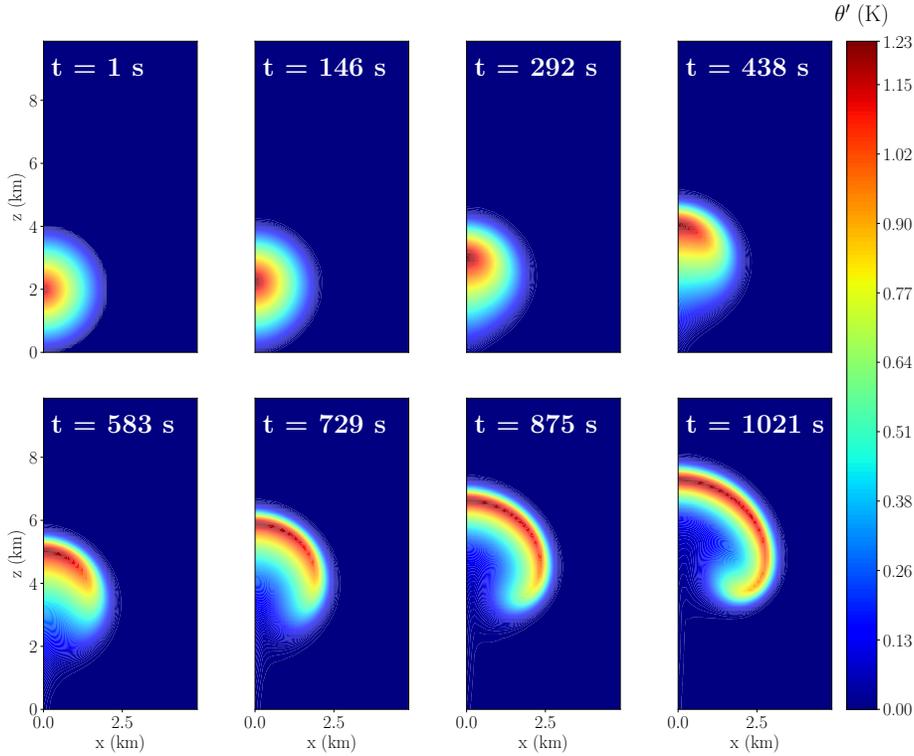}
	\caption{Evolution of the potential temperature perturbation field $\theta'$ for the rising thermal bubble test case.}
	\label{fig:thermal-bubble-evolution}
\end{figure}

We impose impenetrable, free-slip boundary conditions on all domain boundaries. The simulation is run for a total time of $t_f = \SI{1021}{\second}$. During this time interval, the warm bubble rises due to buoyancy forces and evolves into a characteristic mushroom-like shape as a result of shear stresses at its edges. See Figure~\ref{fig:thermal-bubble-evolution}.

For the numerical simulation, we use the artificial viscosity $\mu_a = \SI{15}{\meter\squared\per\second}$ and Prandtl number $\Pr = 1$, as suggested in \cite{Ahmad2006}. These values are chosen ad hoc to stabilize the numerical simulation and ensure the resolution of the main flow features without excessive diffusion.

We generate the reference solution $\theta'_{\text{ref}}(t)$ using $\Delta t_{\text{ref}} = \SI{0.005}{\second}$. For training, we employ coarser time steps $\Delta t_{\text{train}} \in \{\SI{0.01}{\second}, \SI{0.1}{\second}, \SI{1}{\second}\}$. The sets of training and test points are defined as:
\begin{equation*}
\begin{aligned}
    \mathcal{T}_{\text{train}} &= \{t^k_{\text{train}} = k \Delta t_{\text{save}} \mid k = 1, ..., N_T\}, \\
    \mathcal{T}_{\text{test}} &= \{t_{\text{test}}^j = j + (0.05 \pm \Delta t_{\text{test}}) \mid j = 1, ..., N_T\},
\end{aligned}
\end{equation*}
where $\Delta t_{\text{save}} = \SI{1}{\second}$ and $N_T = 1021$.

Table \ref{tab:disc_error} reports
the discretization error $\overline{E}_{\text{disc}}(\Delta t_{\text{train}})$ for each $\Delta t_{\text{train}}$, along with execution times and time ratios relative to the reference case.
The reference case, with $\Delta t = \SI{0.005}{\second}$, has an execution time of \SI{659.76}{\second}. As expected, when going from $\Delta t = \SI{1}{\second}$ to $\Delta t = \SI{0.1}{\second}$, the discretization error decreases. However, the discretization error remains of the same order of magnitude when further reducing the time step to $\Delta t = \SI{0.01}{\second}$, suggesting that the error due to time discretization becomes negligible with respect to the spatial discretization error.
For this reason, we will show only the results for $\Delta t = \SI{0.1}{\second}$ and $\Delta t = \SI{1}{\second}$ in the subsequent analysis.
Table \ref{tab:disc_error} also highlights the computational efficiency gained by using larger time steps. For instance, the simulation with $\Delta t = \SI{1}{\second}$ runs approximately 46 times faster than the reference case, while still maintaining a reasonable level of accuracy.

\begin{table}[ht]
    \centering
    \begin{tabular}{|c|c|c|c|}
        \hline
        $\Delta t_{\text{train}}$ & $\overline{E}_{\text{disc}}(\Delta t_{\text{train}})$ & Execution Time & Time Ratio \\
        \hline
        \SI{1}{\second}  & 3.88e-02 & \SI{14.18}{\second} & 46.53 \\
        \SI{0.1}{\second}  & 8.53e-03 & \SI{45.54}{\second} & 14.81 \\
        \SI{0.01}{\second} & 9.63e-03 & \SI{332.19}{\second} & 1.99 \\
        \hline
    \end{tabular}
    \caption{Discretization error $\overline{E}_{\text{disc}}(\Delta t_{\text{train}})$, execution time, and time ratio relative to the reference case (with $\Delta t = \SI{0.005}{\second}$) for different training time steps.}
    \label{tab:disc_error}
\end{table}

\begin{figure}[ht]
    \centering
    \includegraphics[width=.85\textwidth]{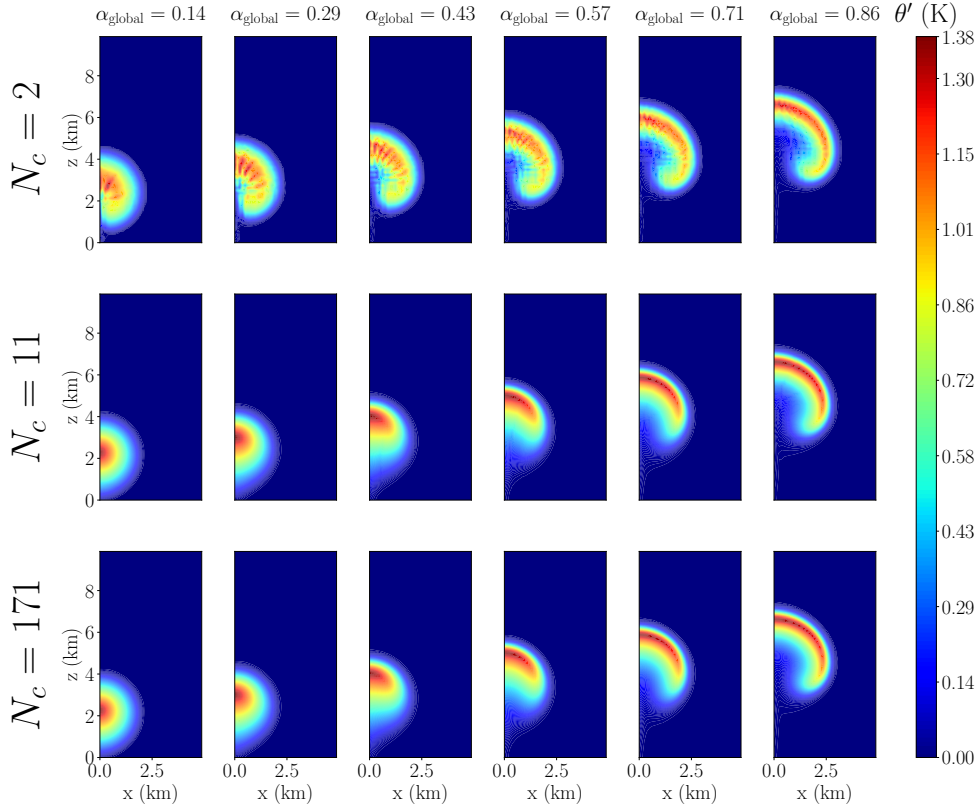}
    \caption{Evolution of synthetic solutions for the rising thermal bubble as the number of checkpoints $N_c$ varies. The images demonstrate how increasing $N_c$ leads to more physically relevant interpolated trajectories.}
    \label{fig:synthetic-evolution}
\end{figure}

In the checkpoint selection process, we vary the number of checkpoints $N_c$ while maintaining a constant total number of synthetic snapshots $N_{\text{tot}} = 1021$. For each $N_c$, we adapt $N_{\text{synth}}$ according to:
\begin{equation*}
N_{\text{synth}} = \left\lfloor\frac{N_{\text{tot}} - N_c}{N_c - 1} \right\rfloor.
\end{equation*}
Figure \ref{fig:synthetic-evolution} illustrates the evolution of synthetic solutions as the number of checkpoints varies for $\Delta t _{\text{train}}= \SI{1}{\second}$. Note that checkpoints act as forced passage points for the interpolated field. For instance, when $N_c = 2$, the interpolant follows a path that causes the supports to move along the diagonal connecting the initial condition to the final solution, resulting in a trajectory
of limited physical relevance. However, as we increase the number of checkpoints, we constrain the solution to follow a progressively more physically meaningful trajectory.

Table \ref{tab:computation_time} shows the percentage breakdown of computation time for different steps of our method, as well as the total computation time, for different numbers of checkpoints.
We see that Steps 2-3, which involves solving the OT problems in \eqref{eq:step2}, carry the majority of the computational cost, especially as the number of checkpoints increases. However, these problems are entirely independent of each other, and thus amenable to complete parallelization\footnote{In our current implementation, we have not fully leveraged this parallelization potential, as the computations were not performed on a distributed cluster.}.

\begin{table}[htbp]
\centering
\begin{tabular}{|l|cccccc|}
\hline
 & \multicolumn{6}{c|}{Number of Checkpoints ($N_c$)} \\
 & 2 & 11 & 31 & 61 & 103 & 171 \\
\hline
Step 2-3 (\%) & 81.74 & 90.79 & 96.80 & 98.28 & 99.09 & 99.44 \\
Step 5 (\%)   & 18.25 & 9.20  & 3.20  & 1.72  & 0.91  & 0.56 \\
Total Time (s) & 10.44 & 30.54 & 91.78 & 178.26 & 309.99 & 505.17 \\
\hline
\end{tabular}
\caption{Computation time breakdown for different numbers of checkpoints ($N_c$) with time step $\Delta t = \SI{1.0}{\second}$.}
\label{tab:computation_time}
\end{table}

\begin{figure}[ht]
    \centering
    \includegraphics[width=0.6\textwidth]{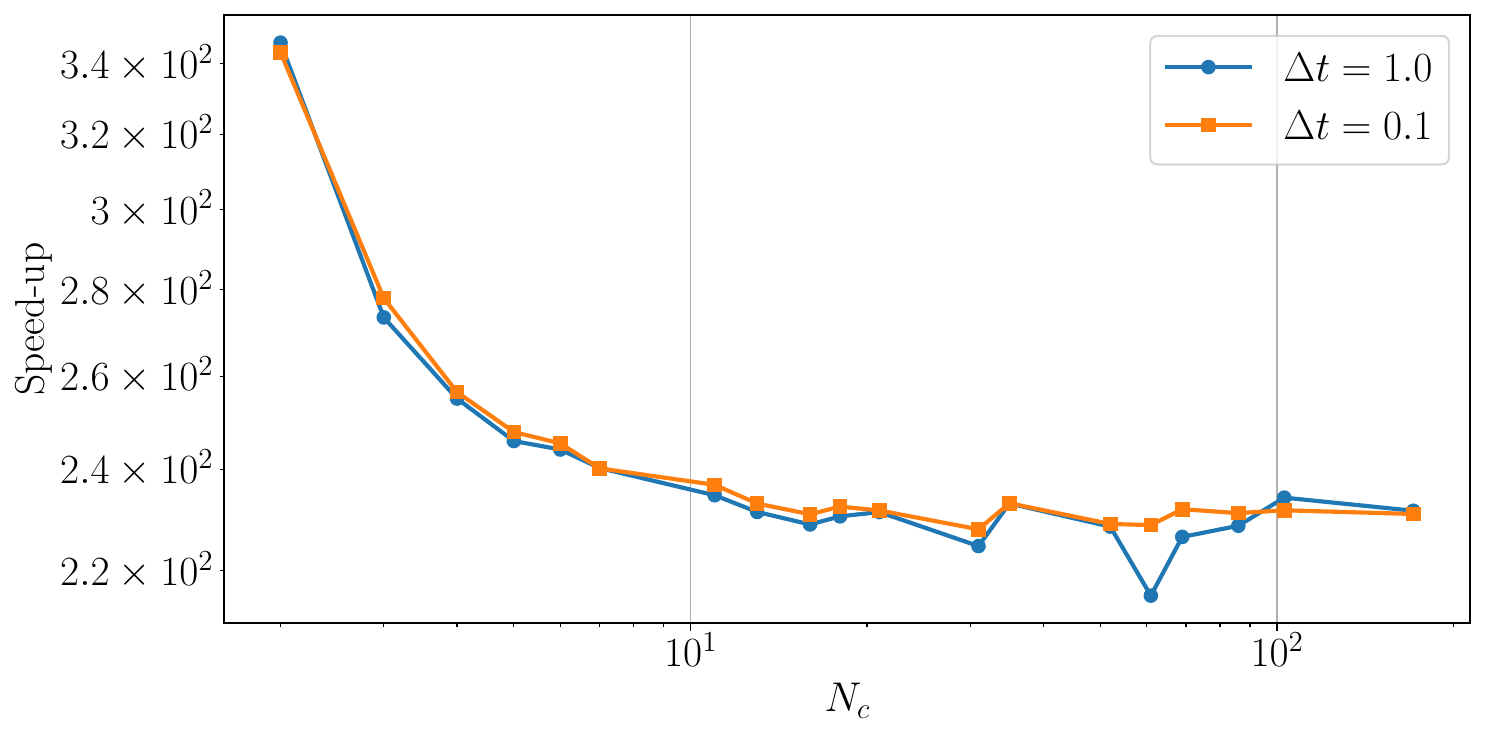}
    \caption{Speed-up  factor as a function of the number of checkpoints $N_c$ for the rising thermal bubble test case.}
    \label{fig:speed_up}
\end{figure}

To assess the computational efficiency of our approach for the \textit{online phase}, we introduce the concept of \textit{speed-up}. This metric is defined as the ratio between the reference FOM time (for $\Delta t_{\text{test}} = \SI{0.005}{\second}$) and the time required for Step 5, which represents the online inference phase of our method. It is important to note that we are primarily comparing the online inference times, as this phase utilizes precomputed OT plans, resulting in a simple matrix-vector multiplication with relatively low computational cost. Figure \ref{fig:speed_up} illustrates how the speed-up varies with the number of checkpoints $N_c$ for different time steps $\Delta t$. The speed-up decreases with increasing $N_c$, but remains advantageous \((\mathcal{O}(10^2))\) for all considered \(N_c\) values.

\begin{figure}[ht]
\centering
\includegraphics[width=1\textwidth]{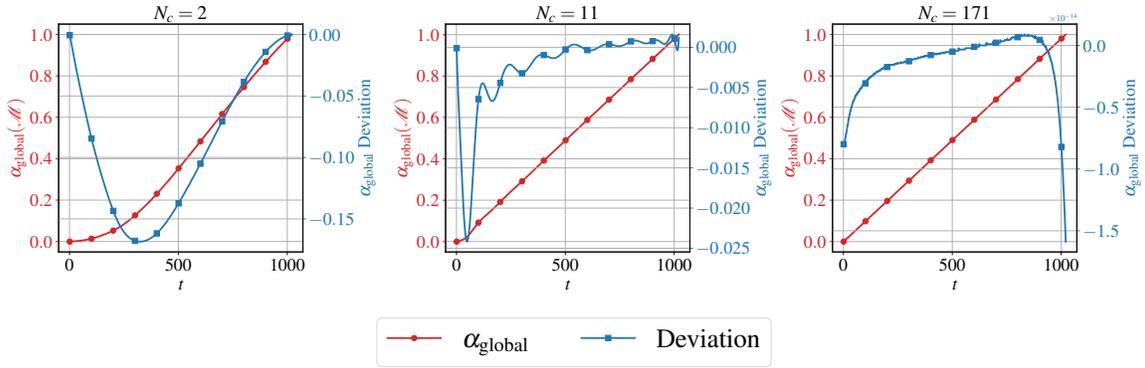}
\caption{Global interpolation parameter $\alpha_{\text{global}}(\mathcal{M})$ in red and its deviation in blue over time for $N_c = 2$ (left), $N_c = 11$ (center), $N_c = 171$ (right).}
\label{fig:alpha_global}
\end{figure}

To better understand the relationship between physical time and the interpolation parameter, in Figure \ref{fig:alpha_global} we examine the behavior of the global interpolation parameter $\alpha_{\text{global}}(\cdot,\mathcal{M})$ as a function of time for different numbers of checkpoints.
We observe that increasing the number of checkpoints causes the $\mathcal{M}$ mapping to progressively approach the linear mapping $\mathcal{L}$, as explained in Remark \ref{remark1}.

\begin{figure}[ht]
\centering
\includegraphics[width=1\textwidth]{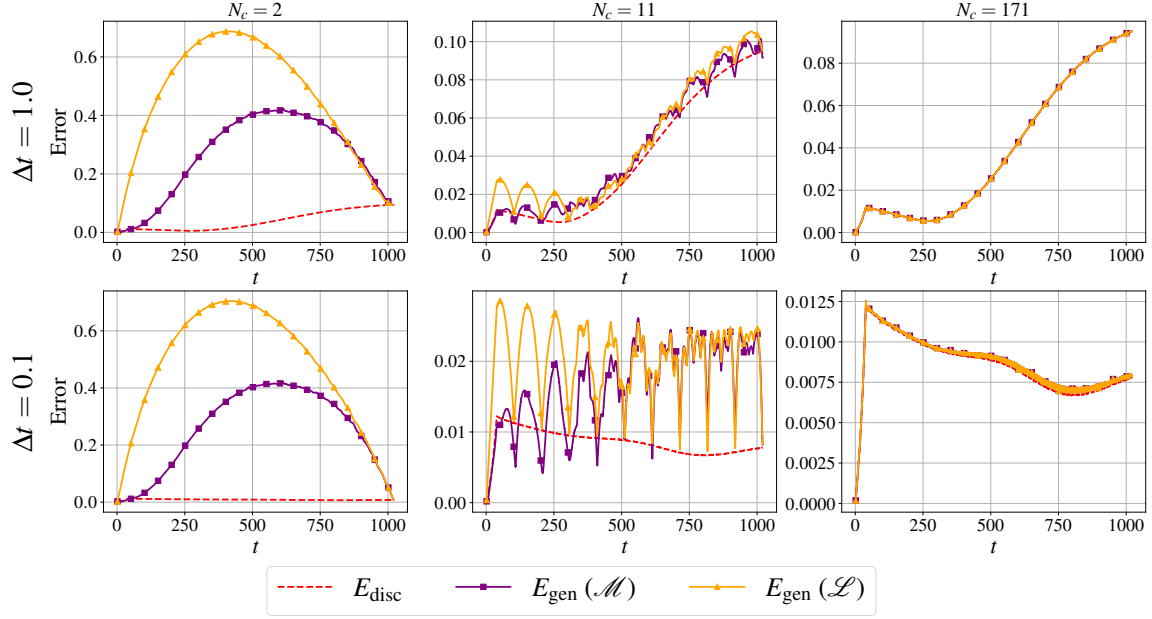}
\caption{Generalization error $E_{\text{gen}}$ over time $t$ for both mappings $\mathcal{M}$ and $\mathcal{L}$ and different values of $N_c$ and $\Delta t$. Note that the scale of the vertical axis varies in each subplot.}
\label{fig:reconstruction_error}
\end{figure}

Figure \ref{fig:reconstruction_error} shows the generalization errors for both mappings $\mathcal{M}$ and
$\mathcal{L}$ and selected numbers of checkpoints and time step sizes. $E_{\text{gen}}$ is made up of two components: one from displacement interpolation and the other
from temporal discretization. The interpolation error decreases with increasing $N_c$, showing oscillatory behavior between checkpoints. Local minima occur near checkpoints, reflecting reduced error at known data points. As $N_c$ increases, this error decreases,
lowering the overall $E_{\text{gen}}$ until the
time discretization error dominates.

\begin{figure}[ht]
\centering
\includegraphics[width=1\textwidth]{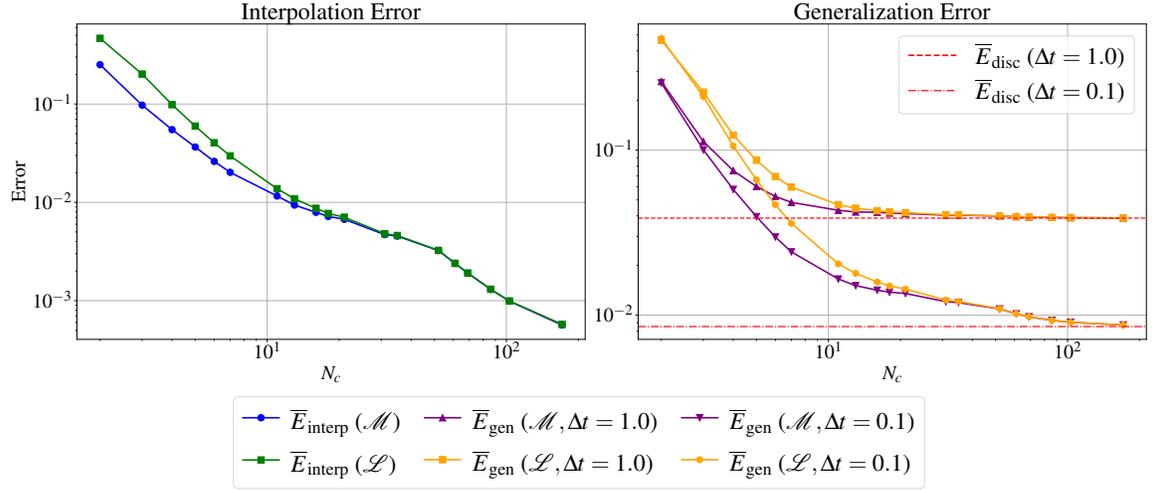}
\caption{Mean interpolation error $\overline{E}_{\text{interp}}$ (left) and mean generalization error $\overline{E}_{\text{gen}}$ (right)
for both mappings $\mathcal{M}$ and $\mathcal{L}$
vs. the number of checkpoints $N_c$ for different time step sizes $\Delta t$.}
\label{fig:mean_error}
\end{figure}

Figure \ref{fig:mean_error} provides a comprehensive view of the mean interpolation and generalization errors as functions of $N_c$ for both mappings $\mathcal{M}$ and $\mathcal{L}$ and different $\Delta t$ values.
Only $E_{\text{interp}}$ for $\Delta t = \SI{1}{\second}$
is presented as it remains consistent across $\Delta t$, indicating robustness to temporal discretization.
This can be attributed to the underlying displacement interpolation capability: since the checkpoints are similar for different $\Delta t$, the interpolant behavior within each interval remains comparable.

Let us summarize the results shown so far.
The mean interpolation error $\overline{E}_{\text{interp}}$ decreases monotonically with increasing $N_c$ for both mappings, with $\mathcal{M}$ consistently outperforming $\mathcal{L}$, especially for lower values of $N_c$, when the deviation from linearity is significant, as shown in Figure \ref{fig:alpha_global}. The mean generalization error $\overline{E}_{\text{gen}}$ exhibits a similar trend, but saturates at the underlying time discretization error $\overline{E}_{\text{disc}}$ for each $\Delta t$.  This saturation highlights an important characteristic of our method: by increasing $N_c$, we can effectively reduce the error associated with the displacement interpolation representation, but we are ultimately bounded by the inherent temporal discretization error of the training data.
Thus, our method is able to achieve optimal accuracy within the constraints of the temporal resolution
of the training data.

\subsubsection*{POD Analysis of Synthetic Datasets}

To further evaluate the effectiveness of our data augmentation strategy, we performed POD on two datasets. The first is the synthetic snapshot matrix $\mathbf{S}_{\text{synth}} \in \mathbb{R}^{N_h \times N_\text{tot}}$, as defined in Sec.~\ref{sec:ot_rom_framework}, which includes both the original checkpoints and the interpolated states. The second is the checkpoint snapshot matrix $\mathbf{S}_{\text{check}} \in \mathbb{R}^{N_h \times N_c}$, containing only the $N_c$ checkpoints selected from the FOM trajectory $\mathcal{S}_h$.
For both cases, we set an energy threshold of 99.99\% to determine the number of POD modes to retain. This high threshold ensures that we capture almost all of the system's energy in our reduced basis.
Let $\mathbf{U}_{\text{synth}}$ and $\mathbf{U}_{\text{check}}$ denote the POD bases obtained from $\mathbf{S}_{\text{synth}}$ and $\mathbf{S}_{\text{check}}$, respectively. We compute
the projection errors for both the training and test datasets onto these bases, which for a given snapshot $\boldsymbol{u}_h(t)$ onto a basis $\mathbf{U}$ is given as\footnote{Here, $\boldsymbol{u}_h(t)$ denotes the vector of degrees of freedom associated with the discretized solution at time $t$.}:
\begin{equation}\label{eq:errproj}
E_{\text{proj}}(\boldsymbol{u}_h(t), \mathbf{U}) = \frac{\|\boldsymbol{u}_h(t) - \mathbf{U}\mathbf{U}^T\boldsymbol{u}_h(t)\|_2}{\|\boldsymbol{u}_h(t)\|_2}.
\end{equation}
We compute $E_{\text{proj}}$ in \eqref{eq:errproj} for both $\mathcal{S}_{\text{train}}$ and $\mathcal{S}_{\text{test}}$, using $\mathbf{U}_{\text{synth}}$ and $\mathbf{U}_{\text{check}}$.

\begin{figure}[htpb]
\centering
\includegraphics[width=1\textwidth]{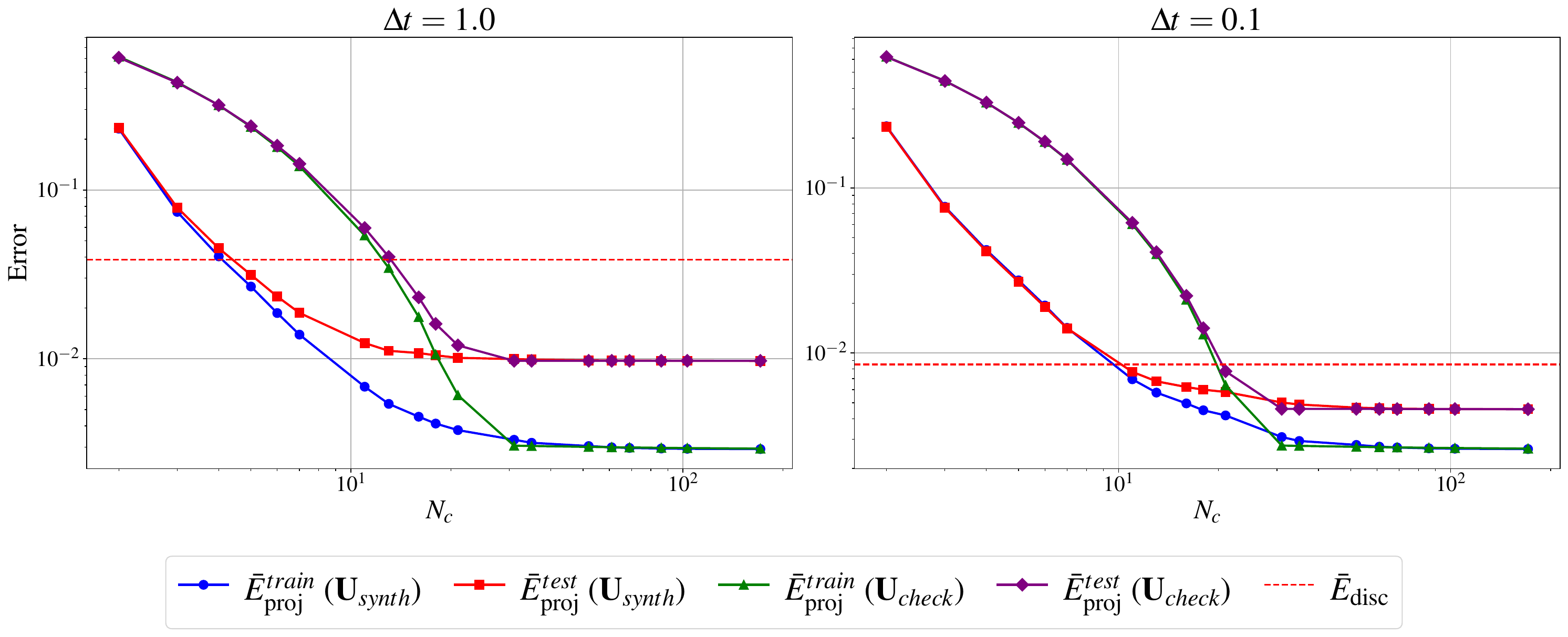}
\caption{Projection errors for POD bases derived from synthetic and checkpoint datasets for $\Delta t = 1$ (left) and $\Delta t = 0.1$ (right).}
\label{fig:pod_projection_errors}
\end{figure}

Figure \ref{fig:pod_projection_errors} presents the projection errors for both $\mathbf{U}_{\text{synth}}$ and $\mathbf{U}_{\text{check}}$ across different values of $N_c$ and $\Delta t$. For $\Delta t = \SI{1}{\second}$, the projection error for $\mathbf{U}_{\text{synth}}$ becomes comparable to $\overline{E}_{\text{disc}}$ with as few as 3 checkpoints, while for $\Delta t = \SI{0.1}{\second}$ it takes 11 checkpoints.
This difference can be attributed to the finer temporal resolution capturing more detailed dynamics
and hence necessitating a larger number of checkpoints to adequately represent the solution manifold,
and the lower $\overline{E}_{\text{disc}}$ threshold for $\Delta t = \SI{0.1}{\second}$ . We see that
$\mathbf{U}_{\text{synth}}$ consistently outperforms $\mathbf{U}_{\text{check}}$, especially for lower $N_c$
values, with this superiority more pronounced for $\Delta t = \SI{1}{\second}$. As $N_c$ increases, errors for both bases converge to similar values.
This happens earlier for $\Delta t = \SI{0.1}{\second}$, suggesting that, for sufficiently high $N_c$, the additional information from synthetic snapshots becomes less critical.

\subsection{Density Current}
\label{subsec:density-current}
\graphicspath{{Images/density_current/}}

The density current test case simulates the evolution of a cold air bubble in a neutrally stratified atmosphere, providing a more challenging benchmark for gravity-driven flows and complex structure formation compared to the thermal bubble case.

\begin{figure}[ht]
    \centering
    \includegraphics[width=\textwidth]{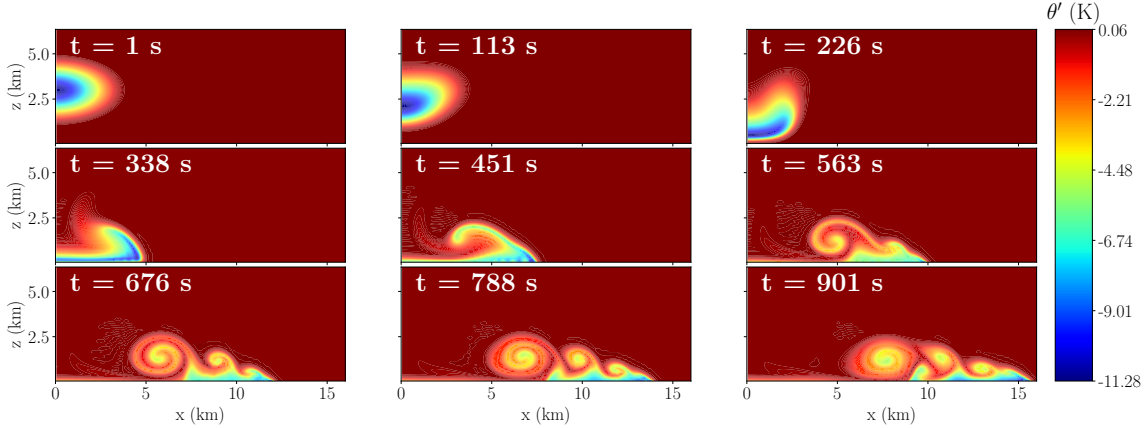}
    \caption{Evolution of the potential temperature perturbation field $\theta'$ for the density current test case.}
    \label{fig:density-current-evolution}
\end{figure}

The computational domain is $\Omega = [\SI{0}{\meter}, \SI{25600}{\meter}] \times [\SI{0}{\meter}, \SI{6400}{\meter}]$, discretized with a uniform structured mesh with mesh size of $h = \Delta x = \Delta z = \SI{100}{\meter}$, resulting in a $256 \times 64$ grid and $N_h = 16384$.
The initial conditions are:

\begin{equation}
\left\{
\begin{aligned}
{\theta}_s &= \SI{300}{\kelvin}, \\[6pt]
\theta_0 &= \begin{cases}
    {\theta}_s - \SI{7.5}{\kelvin}\left[1 + \cos(\pi r)\right] & \text{if } r \leq 1 \\
    {\theta}_s & \text{otherwise}
\end{cases} ,\\[6pt]
\rho_0 &= \frac{p_g}{R\theta_0}\left(\frac{p}{p_g}\right)^{c_v/c_p} \text{ with } p = p_g\left(1 - \frac{gz}{c_p\theta_0}\right)^{c_p/R}, \\[6pt]
\u_0 &= \mathbf{0}, \\[6pt]
\end{aligned}
\right.
\end{equation}
where $r = \sqrt{\left(\frac{x-x_c}{x_r}\right)^2 + \left(\frac{z-z_c}{z_r}\right)^2}$, with $(x_r, z_r) = (\SI{4000}{\meter}, \SI{2000}{\meter})$ and $(x_c, z_c) = (\SI{0}{\meter}, \SI{3000}{\meter})$.
We consider impenetrable, free-slip boundary conditions and run the simulation for $t_f = \SI{901}{\second}$. Figure \ref{fig:density-current-evolution} shows the evolution of the potential temperature perturbation field.

Following \cite{Straka1993}, we set $\mu_a = \SI{75}{\meter\squared\per\second}$ and $\Pr = 1$. The reference solution uses $\Delta t_{\text{ref}} = \SI{0.005}{\second}$, while training solutions use $\Delta t_{\text{train}} \in \{\SI{0.01}{\second}, \SI{0.1}{\second}\}$. Unlike the thermal bubble case, the simulation with $\Delta t = \SI{1}{\second}$ did not converge, indicating higher sensitivity to temporal discretization.

Table \ref{tab:density_current_error} reports the discretization error, execution time, and time ratio for different training time steps. The reference FOM execution time on the test is \SI{2795.76}{\second}.

\begin{table}[ht]
    \centering
    \begin{tabular}{|c|c|c|c|}
        \hline
        $\Delta t_{\text{train}}$ & $\overline{E}_{\text{disc}}(\Delta t_{\text{train}})$ & Execution Time & Time Ratio \\
        \hline
        \SI{0.1}{\second}  & 4.40e-02 & \SI{213.04}{\second} & 13.12 \\
        \SI{0.01}{\second} & 4.76e-03 & \SI{1515.11}{\second} & 1.85 \\
        \hline
    \end{tabular}
    \caption{Discretization error $\overline{E}_{\text{disc}}(\Delta t_{\text{train}})$, execution time, and time ratio relative to the reference case (with $\Delta t = \SI{0.005}{\second}$) for different training time steps.}
    \label{tab:density_current_error}
\end{table}

\begin{figure}[ht]
    \centering
    \includegraphics[width=\textwidth]{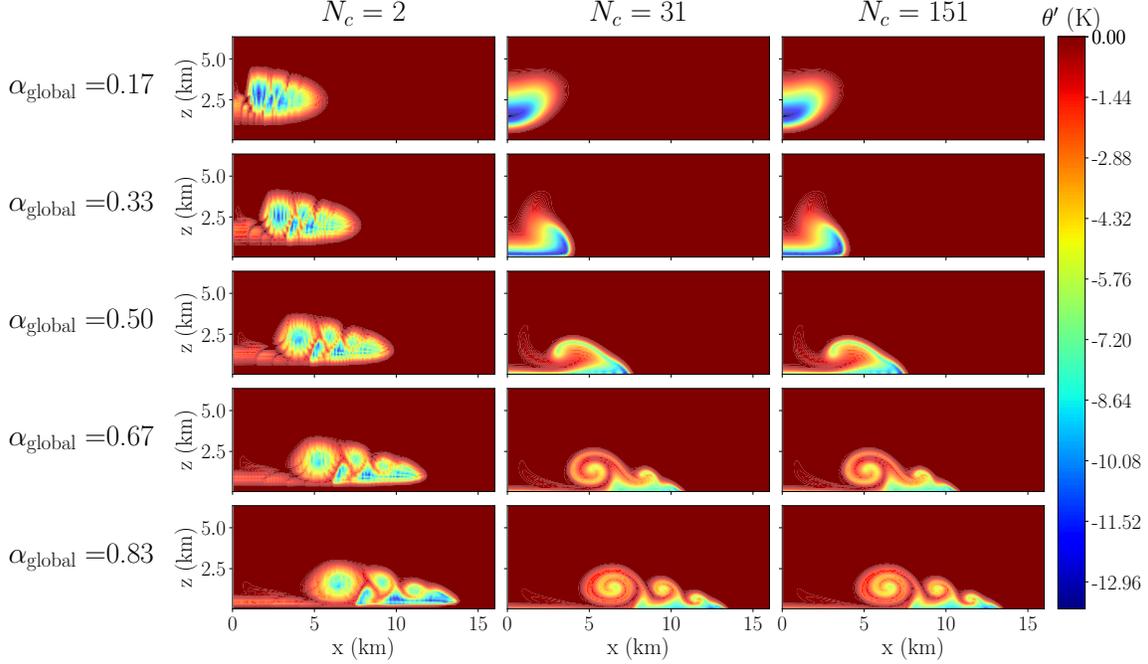}
    \caption{Evolution of synthetic solutions for the density current as the number of checkpoints $N_c$ varies.}
    \label{fig:density-current-barycenter}
\end{figure}

Figure \ref{fig:density-current-barycenter} illustrates the evolution of synthetic solutions
for different values of $N_c$ with $\Delta t_{\text{train}} = \SI{0.1}{\second}$.
The density current evolution exhibits a combination of two distinct motion phases: a mostly vertical motion until approximately $t = \SI{280}{\second}$, followed by horizontal propagation at the ground. This dual-phase motion makes the choice of checkpoints particularly crucial for accurate interpolation, and again we see how increasing $N_c$ leads to more physically relevant interpolated trajectories.

\begin{figure}[ht]
\centering
\includegraphics[width=1\textwidth]{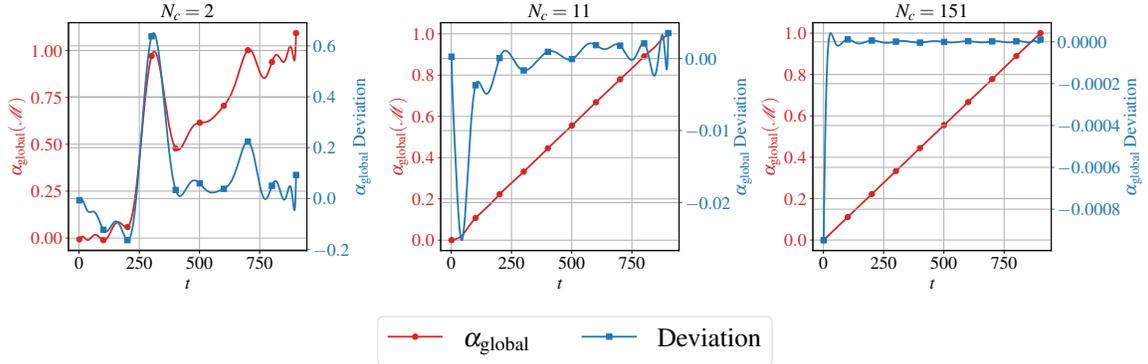}
\caption{Global interpolation parameter $\alpha_{\text{global}}(\mathcal{M})$ in red and its deviation in blue over time for $N_c = 2$ (left), $N_c = 11$ (center), $N_c = 151$ (right).}
\label{fig:density-current-alpha-global}
\end{figure}

Figure \ref{fig:density-current-alpha-global} shows the evolution in time of $\alpha_{\text{global}}(\cdot,\mathcal{M})$. We observe more a pronounced nonlinearity than the thermal bubble case, especially for lower values of $N_c$. This reflects the significant variations in the solution's rate of change throughout the simulation.

\begin{figure}[ht]
\centering
\includegraphics[width=1\textwidth]{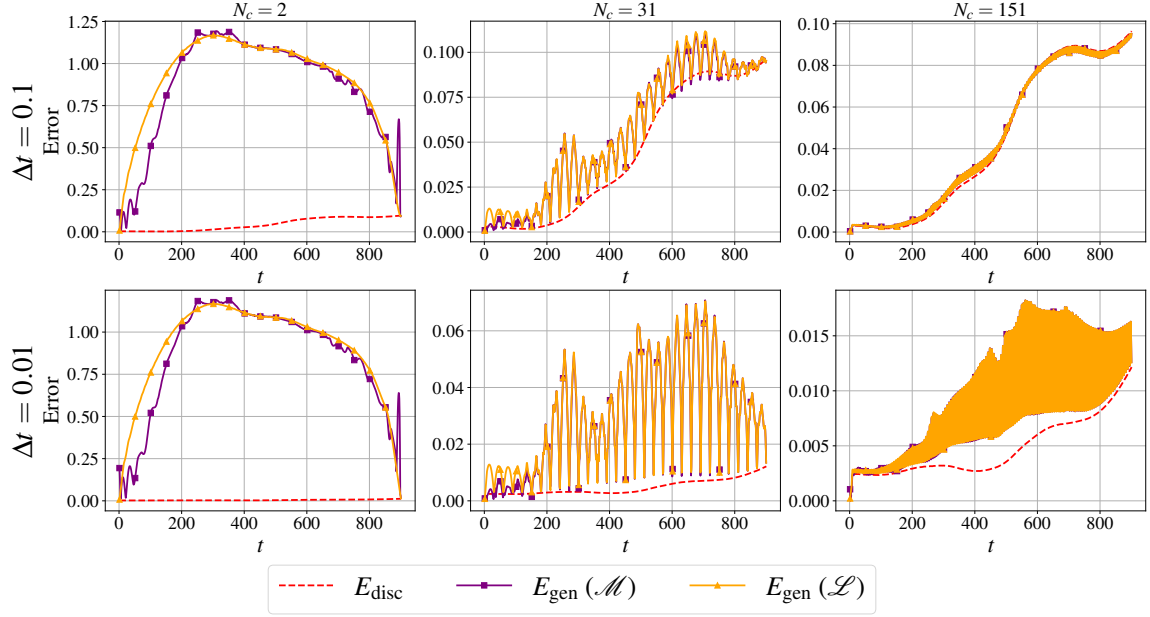}
\caption{Generalization error $E_{\text{gen}}$ over time $t$ for both mappings $\mathcal{M}$ and $\mathcal{L}$ and different values of $N_c$ and $\Delta t$. Note that the scale of the vertical
axis varies in each subplot.}
\label{fig:density-current-generalization}
\end{figure}

Figure \ref{fig:density-current-generalization} presents the generalization error $E_{\text{gen}}$ for both mappings $\mathcal{M}$ and $\mathcal{L}$ and selected values of $N_c$ and $\Delta t$. Compared to the thermal bubble case, we observe higher overall error levels and more pronounced oscillations, particularly for $N_c = 2$. The difference in error magnitude between $\Delta t = \SI{0.1}{\second}$ and $\Delta t = \SI{0.01}{\second}$ is more significant, highlighting the density current's greater sensitivity to temporal discretization.
Interestingly, the errors between mapping $\mathcal{M}$ and mapping $\mathcal{L}$
are of the same order of magnitude, with only a small advantage for $\mathcal{M}$ for small $N_c$ values. This suggests that the regression task for $\alpha_{\text{global}}$ introduces an error that partially offsets the potential accuracy gain of the nonlinear mapping $\mathcal{M}$.

\begin{figure}[ht]
\centering
\includegraphics[width=1\textwidth]{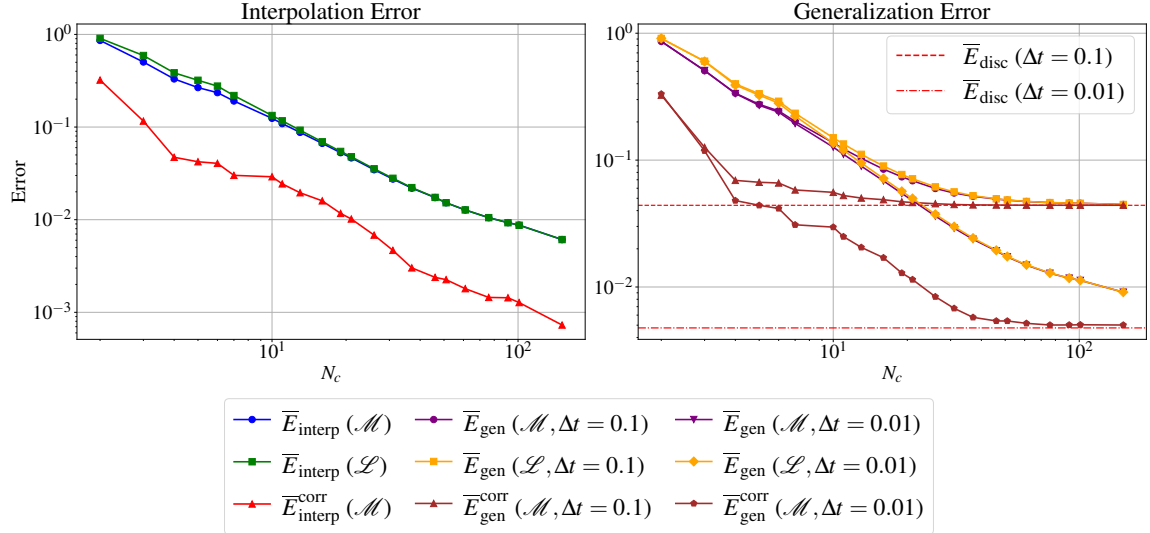}
\caption{Mean interpolation error $\overline{E}_{\text{interp}}$ (left) and mean generalization error $\overline{E}_{\text{gen}}$ (right)
for both mappings $\mathcal{M}$ and $\mathcal{L}$
vs. the number of checkpoints $N_c$ for different time step sizes $\Delta t$.
The goal is to show the effect of the POD-based correction on errors $\overline{E}_{\text{interp}}$ and $\overline{E}_{\text{gen}}$.}
\label{fig:pod-correction}
\end{figure}

To further improve prediction accuracy, we implemented the POD-based correction described in Sec.~\ref{sec:ot_rom_framework}, at point (6). We collected residuals from the barycentric representation for the $\mathcal{S}_{\text{train}}$ dataset and performed POD with an energy threshold of 99.99\%. Figure \ref{fig:pod-correction} shows the effect of this correction on the mean interpolation and generalization error.
We see that the POD-based correction consistently reduces the error across different
values of $N_c$ and $\Delta t$, with the improvement being particularly significant
when using fewer checkpoints. This correction is especially valuable for capturing fine-scale structures like Kelvin-Helmholtz instability along the upper boundary of the cold front, which are challenging to represent accurately through OT-based interpolation alone.

The effectiveness of the POD correction varies with $N_c$ and $\Delta t$. For larger time steps and fewer checkpoints, when the OT-based interpolation provides a coarser approximation, the POD correction plays a crucial role.
As $N_c$ increases or $\Delta t$ decreases, i.e., as the OT-based interpolation itself becomes more accurate, the relative impact of the POD correction diminishes.
These results demonstrate the complementary nature of OT-based interpolation and POD correction in our framework. The OT approach captures large-scale transport phenomena, while the POD correction refines the prediction by accounting for finer-scale dynamics and systematic biases. This combination allows for accurate reduced-order modeling across various temporal resolutions and checkpoint densities, making it particularly suitable for complex atmospheric flows like the density current.

\section{Conclusions and Perspectives}
\label{sec:conclusions}

We presented a novel approach for reduced-order modeling of complex nonlinear dynamical systems, leveraging displacement interpolation in optimal transport theory.
Our method addresses two main challenges that are common to various scientific and engineering applications: capturing nonlinear dynamics and dealing with limited data.
By combining OT-based interpolation with POD with regression, we obtained improved accuracy and efficiency in predicting complex phenomena, as evidenced by our results on challenging benchmark problems.

The proposed framework offers several advantages over traditional ROM techniques:
i) it better captures nonlinear dynamics thanks to the OT-based interpolation;
ii) it includes an efficient data augmentation strategy that allows for improved training of deep learning-based ROMs;
iii) it improves temporal resolution through synthetic snapshot generation; iv) it is flexible in handling
nonlinear systems with limited data.

While our current results are promising, there are several avenues for future research and potential extensions of this work. One potential direction is the development of an adaptive checkpoint selection strategy to reduce the number of required checkpoints while maintaining or even improving accuracy. One could also explore advanced interpolation and extrapolation schemes to enhance the method.
Another perspective is to perform the interpolation with respect to a carefully chosen reference configuration, rather than between consecutive snapshots.
Finally, the integration of our method with real-time data assimilation techniques is another interesting research direction. By continuously updating the ROM with new observational data, one could improve long-term prediction accuracy and adapt to changing system conditions.

\section*{Acknowledgments}
\textbf{MK}, \textbf{FP} and \textbf{GR} acknowledge the support provided by the European Union - NextGenerationEU, in the framework of the iNEST - Interconnected Nord-Est Innovation Ecosystem (iNEST ECS00000043 – CUP G93C22000610007) consortium and its CC5 Young Researchers initiative.
The authors also like to acknowledge INdAM-GNSC for its support.

\bibliographystyle{habbrv}
\bibliography{references}
\end{document}